\documentclass[11pt]{article}
\usepackage{authblk}

\usepackage{geometry}
\usepackage{amsmath,amssymb,mathrsfs}
\usepackage{color}
\usepackage{graphicx}
\usepackage{subfigure}
\usepackage{soul}
\usepackage{bbding}
\usepackage{lineno,hyperref}
\usepackage{multirow, booktabs}
\modulolinenumbers[5]

\def\bp{\mathbf{p}}
\def\bx{\mathbf{x}}

\def\grad{\nabla}

\def\mS{\mathbb{S}}

\def\mR{\mathbb{R}}

\title{A shallow physics-informed neural network for solving partial differential equations on surfaces}

\author[1,3]{Wei-Fan Hu}
\author[2]{Yi-Jun Shih}
\author[2,3]{Te-Sheng Lin}
\author[2]{Ming-Chih Lai}

\affil[1]{Department of Mathematics, National Central University, Taoyuan 32001, Taiwan}
\affil[2]{Department of Applied Mathematics, National Yang Ming Chiao Tung University, Hsinchu 30010, Taiwan}
\affil[3]{National Center for Theoretical Sciences, National Taiwan University, Taipei 10617, Taiwan}

\begin{document}

\maketitle

\begin{abstract}
In this paper, we introduce a shallow (one-hidden-layer) physics-informed neural network for solving partial differential equations  on static and evolving surfaces.  For the static surface case,  with the aid of level set function, the surface normal and mean curvature used in the surface differential expressions can be computed easily. So instead of imposing the normal extension constraints used in literature, we write the surface differential operators in the form of traditional Cartesian differential operators and use them in the loss function directly.
We perform a series of performance study for the present methodology by solving Laplace-Beltrami equation and surface diffusion equation on complex static surfaces. With just a moderate number of neurons used in the hidden layer, we are able to attain satisfactory prediction results.
Then we extend the present methodology to solve the advection-diffusion equation on an evolving surface with given velocity. To track the surface, we additionally introduce a prescribed hidden layer to enforce the topological structure of the surface and use the network to learn the homeomorphism between the surface and the prescribed topology. The proposed network structure is designed to track the surface and solve the equation simultaneously. Again, the numerical results show comparable accuracy as the static cases. As an application, we simulate the surfactant transport on the droplet surface under shear flow and obtain some physically plausible results.
\end{abstract}

%%%%%%%%%%%%%%%%%%%%%%%%%%%%%%%%%%%%%%%%%
\section{Introduction}
%%%%%%%%%%%%%%%%%%%%%%%%%%%%%%%%%%%%%%%%%

Surface partial differential equations (PDEs) arise in a wide variety of scientific and engineering applications. These equations are formulated in terms of differential operators acting on curved surfaces. Mathematically, they are examples of partial differential equations on manifolds.
Problems of interest include, for instance, modeling of surface-active agents~\cite{HLM18}, deforming vesicles~\cite{AMMV05}, cell motility and chemotaxis~\cite{ESV12}, modeling of biomembranes~\cite{ES10}, restoring a damaged pattern on surfaces~\cite{BBS01}, image processing~\cite{TQZY05}, and computer graphics~\cite{AW13}, etc.

Solving PDEs on surfaces is certainly of major interest among the scientific computing community. The fundamental difficulty comes from the numerical approximation of differential operators along a surface. This long-standing problem has been explored by many researchers for decades. For instance, surface finite element method~\cite{DE07, DE13} is particularly designed for PDEs on discretized triangular surfaces; while generating those triangulation nodes can be time-consuming and the accuracy of the method is significantly affected by the quality of triangulations. Using parametric representation is another natural idea~\cite{O18, GA18}, in which the solution along a smooth surface can be represented via spherical harmonics expansion. However, it may suffer from the intrinsic singularities that are built into the PDE formulation (e.g., poles in spherical coordinates) or in the boundary integral kernel involving Green's function; thus, it needs careful treatments near the singularities. A mesh-free approach called radial basis functions (RBFs) method~\cite{ARM18, WK20} works by first representing the solution by a linear combination of RBFs, and then substituting the approximation at some chosen points into the differential equation directly. As a result, a dense linear system of coefficients must be solved which is likely to be ill-conditioned. In order to obtain a well-conditioned resultant matrix and achieve desired accuracy, it often requires to artificially tune the shape parameter appearing in a certain type of radial basis functions. However, finding such a parameter to reach optimal results remains an issue in the usage of radial basis functions.

On the other hand, embedding techniques solve the PDE in a small band in the vicinity of the surface, examples including level set method~\cite{BCO01}, closest point method~\cite{RM08, PR16, PLPR19}, or grid based particle method~\cite{LZ09, LLZ11}.
The underlying surface PDE is alternatively represented in Eulerian coordinates and thus surface derivatives are replaced by projections of derivatives in the embedding Euclidean space. In such a way, the difficulties such as parameterized or triangulated surfaces can be avoided. Although these methods have the feature of being geometrically flexible, computations in this domain may require imposing suitable conditions at the band's boundary which remains unclear in practice.  Besides, these methods require finding surface projection points of regular Cartesian grids. This task needs further computational efforts and can be troublesome when a highly oscillatory surface is considered.

As far as we know, only a few works exist using machine learning for solving PDEs on surfaces. Following the same spirit as those in embedding techniques, Fang et al.~\cite{FZ19, FZY21} adopts the physics-informed neural networks (PINNs)~\cite{raissi2019} framework to solve the Laplace-Beltrami equation (stationary) and diffusion equation (time-dependent) on static surfaces. The neural network solution is constrained to have zero normal derivatives at given training points along the surface. This restriction leads to an approximate normal extension solution in a narrow band of the surface so the Laplace-Beltrami operator is replaced by the conventional Laplace operator. Hence, the PINNs loss function penalizes the equation residual, and the zero normal first- and second-order derivatives. In such a way, the PDE information only comes from the training points given on the surface so it is significantly different from those embedding techniques. Furthermore, it is completely mesh-free that differs from the aforementioned grid-based embedding methods.

\textcolor{black}{While the loss function seems to be legitimate, the numerical experiments shown in these works~\cite{FZ19, FZY21} have the relative $L^2$ errors more than $1\%$ even deep neural networks are used. So instead of using the above loss function, in this paper, we write the surface differential operators in the form of traditional Cartesian differential operators and use them in the loss function directly. Thus, we can encode the entire embedding PDE without imposing the normal extension constraint. Besides, we adopt a completely shallow (one-hidden-layer) network under PINNs framework so it is easy to implement and train.} \textcolor{black}{As discussed in \cite{Cybenko1989, Hor91, Mhaskar96}, a shallow neural network can theoretically approximate smooth functions and their derivatives accurately. This is the legitimate reason why it can help to solve PDEs in the first place.  The shallow PINNs (or Ritz) method with augmented inputs have been proven very effective for solving elliptic interface problems with jump discontinuities across the interface, see the authors' recent papers in \cite{HLL22, LCLHL22, TLHL22}. }

\textcolor{black}{Until very recently, Tang et al.~\cite{TFR22} proposed a methodology that exactly shares the same spirit as ours, i.e., embedding the solution into Eulerian coordinates and expressing surface differential operators by conventional Cartesian ones in the PINN loss.  However, their numerical experiments (for solving the stationary advection-diffusion equation) adopt the deep network architecture (depth $ = 4$ and width $ = 50$), resulting in numerous parameters to be learned, to reach relative $L^2$ error of magnitude $O(10^{-4})$ for some smooth solutions. By contrast,
we simply use shallow network structures with $60$  neurons to solve stationary and time-dependent PDEs on surfaces that attain satisfactory prediction results with errors $O(10^{-6})$ for the Laplace-Beltrami equation and $O(10^{-5})$ for the surface diffusion equation. Furthermore, we have extended our methodology to
the evolving surface case while the paper in \cite{TFR22} is only a focus on the static surface.}
%To solve PDEs on evolving surfaces, we propose a neural network homeomorphism mapping to track evolving surfaces and use this mapping to evaluate surface geometric quantities.}

The rest of the paper is organized as follows. In Section~\ref{Sec:stationary_PDE}, we first describe a shallow PINNs model to solve stationary PDEs (by taking Laplace-Beltrami equation as an example) on a static surface and perform a series of numerical accuracy tests and comparisons. Then we develop the network solver to solve time-dependent PDEs (by taking diffusion equation with a source term as an example) and also demonstrate its capability for finding solutions on complex surfaces in Section~\ref{Sec:time_PDE}.  In Section~\ref{Sec:evolving_PDE}, we have extended the present methodology to solve the  advection-diffusion equation on 2D evolving surface in $\mR^3$. Some concluding remarks and future works are given in Section~\ref{Sec:conclusion}.

%%%%%%%%%%%%%%%%%%%%%%%%%%%%%%%%%%%%%%%%%
\section{Stationary PDEs on surfaces}\label{Sec:stationary_PDE}
%%%%%%%%%%%%%%%%%%%%%%%%%%%%%%%%%%%%%%%%%

%Let us start by solving stationary partial differential equations on surfaces using neural network approximation first.
Denoting a regular (or smooth) surface by $\Gamma$ embedded in Euclidean space $\mathbb{R}^3$, the considered PDEs take the general form
\begin{align}\label{Eq:PDE1}
\mathcal{L}(u) = f \quad \mbox{ on } \Gamma,
\end{align}
where, for simplicity, $\Gamma$ is assumed to be a closed surface. The operator $\mathcal{L}$ may consist of common differential terms related to the surface geometry, such as surface gradient $\nabla_s u$, surface divergence $\nabla_s\cdot\mathbf{v}$ for some vector field $\mathbf{v}$, or Laplace-Beltrami (or surface Laplace) operator $\Delta_s u$. With suitable surface parametric representation, these differential operators can be evaluated via first and second fundamental forms of differential geometry \cite{Walker15}. For the case that the surface is not closed, some suitable boundary conditions along $\partial\Gamma$ must be given. Nevertheless, the boundary condition does not change the main ingredient of the present methodology (see next subsection).
%In this section, we only focus on time-independent PDEs, while time-dependent problems will be discussed later in Sec.~\ref{Sec:time_PDE}.

As aforementioned, the differential operator $\mathcal{L}$ can be computed using a local parametrization, say $u = u(\theta,\phi)$, where $\theta$ and $\phi$ being surface parameters. However, numerical differentiations of $\mathcal{L}(u)$ using surface parametrization might cause severe numerical instability. For instance, if the considered surface geometry is complicated (a stationary highly oscillatory surface case), or the discretized Lagrangian points are clustered on certain parts of the surface (an time-evolving surface case) can lead to inaccurate computations on derivatives~\cite{VRBZ11}. For the latter case, re-parametrization technique is often required to redistribute those markers on the surface to maintain the numerical accuracy and stability~\cite{VRBZ11, SHL18}.

Our goal is to develop a robust \emph{mesh-free} numerical method for solving PDEs~(\ref{Eq:PDE1}) based on neural network learning technique. To compute $\mathcal{L}(u)$, rather than using surface parametrization, here, we adopt an alternative way using conventional differential operators.
To this end, the solution $u$ defined on the surface is now regarded as an embedded function, $u(x,y,z)$, in the Eulerian space that satisfies $u(x,y,z)=u(\theta,\phi)$ when $(x,y,z)\in\Gamma$. Despite this assumption results in the solution with one dimension higher in variable space (surface coordinates $(\theta,\phi)$ to Cartesian coordinates $(x,y,z)$), the surface differential terms in $\mathcal{L}$ can be rewritten via conventional differential operators in Eulerian coordinates. More precisely, at a given point $\mathbf{x} = (x,y,z)\in\Gamma$, we have
\begin{align}\label{Eq:diff_operator}
\begin{split}
& \nabla_s u = (I-\mathbf{n}\mathbf{n}^T)\nabla u, \\
& \nabla_s\cdot\mathbf{v} = \nabla\cdot\mathbf{v} - \mathbf{n}^T(\nabla\mathbf{v})\mathbf{n}, \\
& \Delta_s u = \Delta u - 2H \partial_nu - \mathbf{n}^T(\nabla^2u)\mathbf{n}.
\end{split}
\end{align}
Here, $\mathbf{n} = \mathbf{n}(\mathbf{x})$ is the unit outward  normal vector on $\Gamma$, $H = H(\mathbf{x})$  is the mean curvature, $\partial_nu = \nabla u \cdot\mathbf{n}$ denotes the normal derivative,  and $\nabla^2u$ is the Hessian matrix of $u$. The derivation of above identities can be found in Appendix. Note that, both normal vector and mean curvature in the above formulas can be directly computed once the level set representation $\psi$ of the surface $\Gamma$ is available. That is, at $\mathbf{x}\in\Gamma$ (so the level set $\psi(\mathbf{x}) = 0$ represents $\Gamma$), the above two geometric quantities can be computed by
\begin{align}\label{Eq:lvset}
\mathbf{n} = \frac{\nabla\psi}{\|\nabla\psi\|} \quad \mbox{and} \quad 2H = \nabla\cdot\mathbf{n} = \frac{\mbox{tr}(\nabla^2\psi) - \mathbf{n}^T(\nabla^2\psi)\mathbf{n}}{\|\nabla\psi\|},
\end{align}
where $\mbox{tr}(\cdot)$ gives the trace of a matrix and $\|\cdot\|$ denotes the standard Euclidean norm.

Throughout the rest of this section, we will only focus on the Laplace-Beltrami equation as
\begin{align}\label{Eq:LB}
 \Delta_s u(\mathbf{x}) = f(\mathbf{x}) \quad \mbox{ on } \Gamma,
\end{align}
where we deliberately put the variable $\mathbf{x}$ to clarify that the differential equation is defined in Eulerian coordinates.
In addition, it is important to mention that, there exists infinitely many embedded functions $u(\mathbf{x})$ whose restriction on $\Gamma$ serves as a solution to Eq.~(\ref{Eq:LB}) (or more generally, Eq.~(\ref{Eq:PDE1})) so that such embedded solutions can be representable in a wide range of neural network approximator thanks to the expressive power of universal approximation theory~\cite{Cybenko1989, Hor91}.

%%%%%%%%%%%%%%%%%%%%%%%%%%%%%%%%%%%%%%%%%
\subsection{Physics-informed learning machinery  using shallow neural network approximation}
%%%%%%%%%%%%%%%%%%%%%%%%%%%%%%%%%%%%%%%%%

With the expressive capabilities of neural networks \cite{Hor91}, we hereby construct a simple feedforward, fully-connected, shallow (one-hidden-layer) neural network approximate solution $u_\mathcal{N}$ as
\begin{align}\label{Eq:NN}
u_\mathcal{N}(\mathbf{x}) = \sum_{j=1}^N \alpha_j \sigma(W_j\mathbf{x}^T+b_j).
\end{align}
Here, $\sigma$ is the activation function, $N$ is the number of employed neurons in that hidden layer. The weights $\alpha_j\in\mathbb{R}$ and $W_j\in\mathbb{R}^{1\times3}$  and the bias $b_j\in\mathbb{R}$  are formed as learnable parameters whose total number is counted as $N_p = 5N$.
Notice that, the output layer in the present network structure is considered to be unbiased so the network output can be concisely written in the form of finite linear combination of activation functions.

Let us describe the methodology of physics-informed learning machinery \cite{raissi2019} for solving  Eq.~(\ref{Eq:LB}) as follows.
With a given training set $\{\mathbf{x}^i = (x^i,y^i,z^i)\in\Gamma\}_{i=1}^M$, the neural net parameters (weights and biases in Eq.~(\ref{Eq:NN})) are learned via minimizing the mean squared error of the differential equation residual
\begin{align*}
\mbox{Loss}(\mathbf{p}) = \frac{1}{M}\sum_{i=1}^M \left[ \Delta_s u_\mathcal{N}(\mathbf{x}^i;\mathbf{p}) - f(\mathbf{x}^i)\right]^2,
\end{align*}
where $\mathbf{p}$ is a vector collecting all training parameters (of dimension $N_p$) .
 Using the third identity in Eq.~(\ref{Eq:diff_operator}), it is
natural to choose the loss function as
\begin{align}\label{Eq:Loss1}
\mbox{Loss}_{\Delta_s}(\mathbf{p}) = \frac{1}{M}\sum_{i=1}^M \left[\Delta u_\mathcal{N}(\mathbf{x}^i) - 2H(\mathbf{x}^i)\partial_nu_\mathcal{N}(\mathbf{x}^i) - \mathbf{n}(\mathbf{x}^i)^T\left(\nabla^2u_\mathcal{N}(\mathbf{x}^i)\right)\mathbf{n}(\mathbf{x}^i) - f(\mathbf{x}^i)\right]^2,
\end{align}
where we have dropped the notation $\mathbf{p}$ in $u_\mathcal{N}$ for succinct purpose. The first- and second-order partial derivatives to $u_\mathcal{N}$ involved in the above loss can be evaluated via auto-differentiation~\cite{BPRS18}, or, derived explicitly through the network expression (\ref{Eq:NN}) thanks to the simplicity of shallow network structure. We remark that the explicit evaluations of partial derivatives can be done more efficiently than the auto-differentiation since the latter one requires multiple runs of backpropagation.

Here, we should point out that the following loss function is used in~\cite{FZY21}
\begin{align}\label{Eq:Loss2}
\mbox{Loss}_{\Delta}(\mathbf{p}) = \frac{1}{M}\sum_{i=1}^M \left[\Delta u_\mathcal{N}(\mathbf{x}^i) - f(\mathbf{x}^i)\right]^2 + \frac{1}{M}\sum_{i=1}^M\left[\partial_nu_\mathcal{N}(\mathbf{x}^i)\right]^2 + \frac{1}{M}\sum_{i=1}^M\left[\mathbf{n}(\mathbf{x}^i)^T\left(\nabla^2u_\mathcal{N}(\mathbf{x}^i)\right)\mathbf{n}(\mathbf{x}^i)\right]^2.
\end{align}
The above loss is inspired by the inequality
\begin{align}\label{Eq:inequality}
|\Delta_s u - f| \leq |\Delta u - f| + |2H||\partial_n u| + |\mathbf{n}^T(\nabla^2 u)\mathbf{n}|
\end{align}
which is also a direct result from the third identity in Eq.~(\ref{Eq:diff_operator}). One can see their idea of designing $\mbox{loss}_\Delta$ in Eq.~(\ref{Eq:Loss2}) is to penalize each term on the righthand side of the above inequality. This will generally result in a normal extension solution (in a very narrow region) since it is attempted to enforce $\partial_n u_\mathcal{N} = 0$ on $\Gamma$. Like our proposed $\mbox{Loss}_{\Delta_s}$ in Eq.~(\ref{Eq:Loss1}), second partial derivatives are still required in Eq.~(\ref{Eq:Loss2}). One favorable feature of $\mbox{Loss}_\Delta$ is to avoid computing the local mean curvatures at training points which will save some computational efforts.
However, since the normal derivative term $\partial_n u$ in the inequality~(\ref{Eq:inequality}) is multiplied by the factor $|2H|$, one can anticipate that the actual differential equation residual using $\mbox{loss}_\Delta$ (\ref{Eq:Loss2}) may become significantly large when an oscillatory surface is considered ($H\gg1$). In next subsection, we will demonstrate that our proposed loss function (\ref{Eq:Loss1}) indeed outperforms the splitting residual loss (\ref{Eq:Loss2}) in the sense of higher predictive accuracy regardless of the surface geometries.
We also point out that since only one hidden layer with moderate number of neurons employed in the present network, the computational complexity and learning workload can be significantly reduced without sacrificing the accuracy.

%%%%%%%%%%%%%%%%%%%%%%%%%%%%%%%%%%%%%%%%%
\subsection{Numerical results}
%%%%%%%%%%%%%%%%%%%%%%%%%%%%%%%%%%%%%%%%%

\begin{figure}[h]
\centering
\includegraphics[scale=0.33]{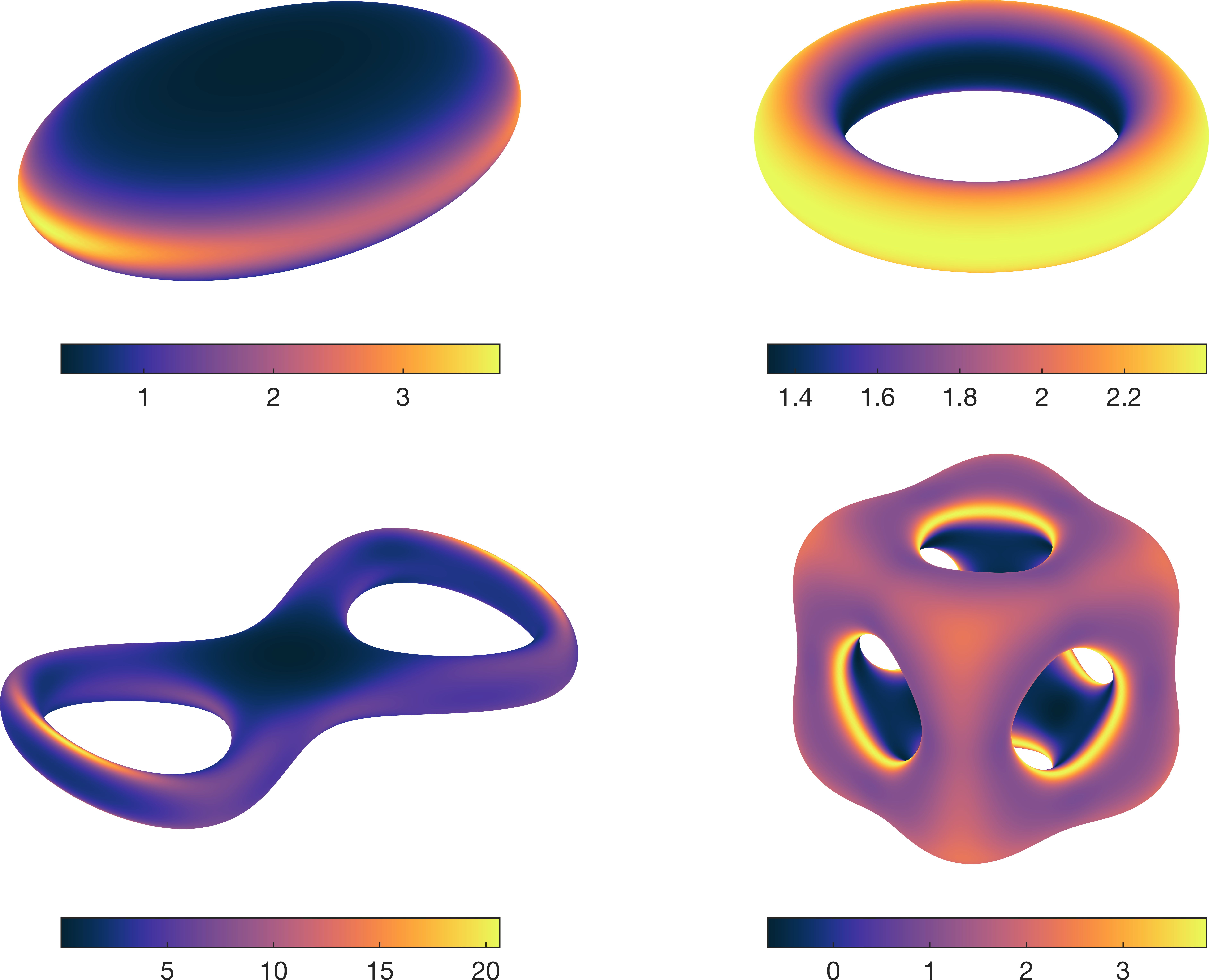}
\caption{Shapes of ellipsoid and torus (top row), genus-2 torus and cheese-like (bottom row). The color code denotes the magnitude of mean curvature $H$.}
\label{Fig:curvature}
\end{figure}

Here we use the established network model to perform a series of numerical tests for Laplace-Beltrami equation. We consider four different geometries of surface $\Gamma$ which can be represented by the zero level set as follows.
\begin{itemize}
\item Ellipsoid : $\psi(x,y,z) = (x/1.5)^2 + (y)^2 + (z/0.5)^2 - 1$
\item Torus: $\psi(x,y,z) = ( \sqrt{x^2+y^2} - 1 )^2 + z^2 - 1/16$
\item Genus-2 torus: $\psi(x,y,z) = \left[ (x+1)x^2(x-1) + y^2 \right]^2 + z^2  - 0.01$
\item Cheese-like surface: $\psi(x,y,z) = (4x^2-1)^2 + (4y^2-1)^2 + (4z^2-1)^2 + 16(x^2+y^2-1)^2 + 16(x^2+z^2-1)^2 + 16(y^2+z^2-1)^2 - 16$
\end{itemize}
As mentioned before, the normal vector and mean curvature used in the computation of Laplace-Beltrami operator $\Delta_s$ can be exactly obtained through symbolic differentiation in Eq.~(\ref{Eq:lvset}). The shapes of these surfaces and corresponding local mean curvatures are shown in Fig.~\ref{Fig:curvature}.

We should note that the solution to Laplace-Beltrami equation is unique up to an arbitrary additive constant, to assess the accuracy of our method, the obtained network solution $u_\mathcal{N}$ is shifted to have the same value of the exact solution at a given point.

Throughout all numerical tests in this paper, we choose sigmoid as the activation function.
We generate a set of collecting points on $\Gamma$ via the usage of \texttt{DistMesh} package developed in~\cite{PS04} wherein level set function related to target surface is required as an input. We then randomly pick training points $\{\mathbf{x}^i\}$ in that point cloud set.
To train the network model, we adopt the Levenberg-Marquardt (LM) method~\cite{Marquardt63} (except the below discussion on the comparison between different popular optimizers) that can effectively find the optimal parameters for losses of mean squared type.
After the training process is finished, we measure the accuracy of the solution
using the test error instead of the training error. That is, we randomly choose $M_{test}$ testing points on $\Gamma$ by computing the relative error in $L^2$ norm as
\begin{align*}
\frac{\|u_\mathcal{N}-u\|_2}{\|u\|_2} = \sqrt{\sum_{i=1}^{M_{test}}\left( u_\mathcal{N}(\mathbf{x}^i) - u(\mathbf{x}^i) \right)^2}\Bigg/\sqrt{\sum_{i=1}^{M_{test}}\left(u(\mathbf{x}^i) \right)^2}.
\end{align*}
For each case, we set $M_{test} = O(10^4)$. And for each test, we repeat the numerical runs for 5 times so the test error reported here is the averaged one.

In the following,  we aim to analyze the performance of our proposed method. We quantify the prediction accuracy through a series of experimental studies, including the comparisons of loss functions and optimizers, and single and double precision computations. We also study the effects on the number of training points and the depth of network architecture.  In the above tests, the ellipsoidal surface is considered, along which the exact solution is chosen as $u(x,y,z) = \sin(x)\cos(y-z)$ so the corresponding right-hand side function $f(x,y,z)$ can be computed directly by substituting $u$ into Eq.~(\ref{Eq:LB}). Furthermore, we also apply the present method to a non-closed surface case (the boundary condition is taken into account) and other more complex surfaces described earlier.

\begin{figure}[h]
\centering
\includegraphics[scale=0.41]{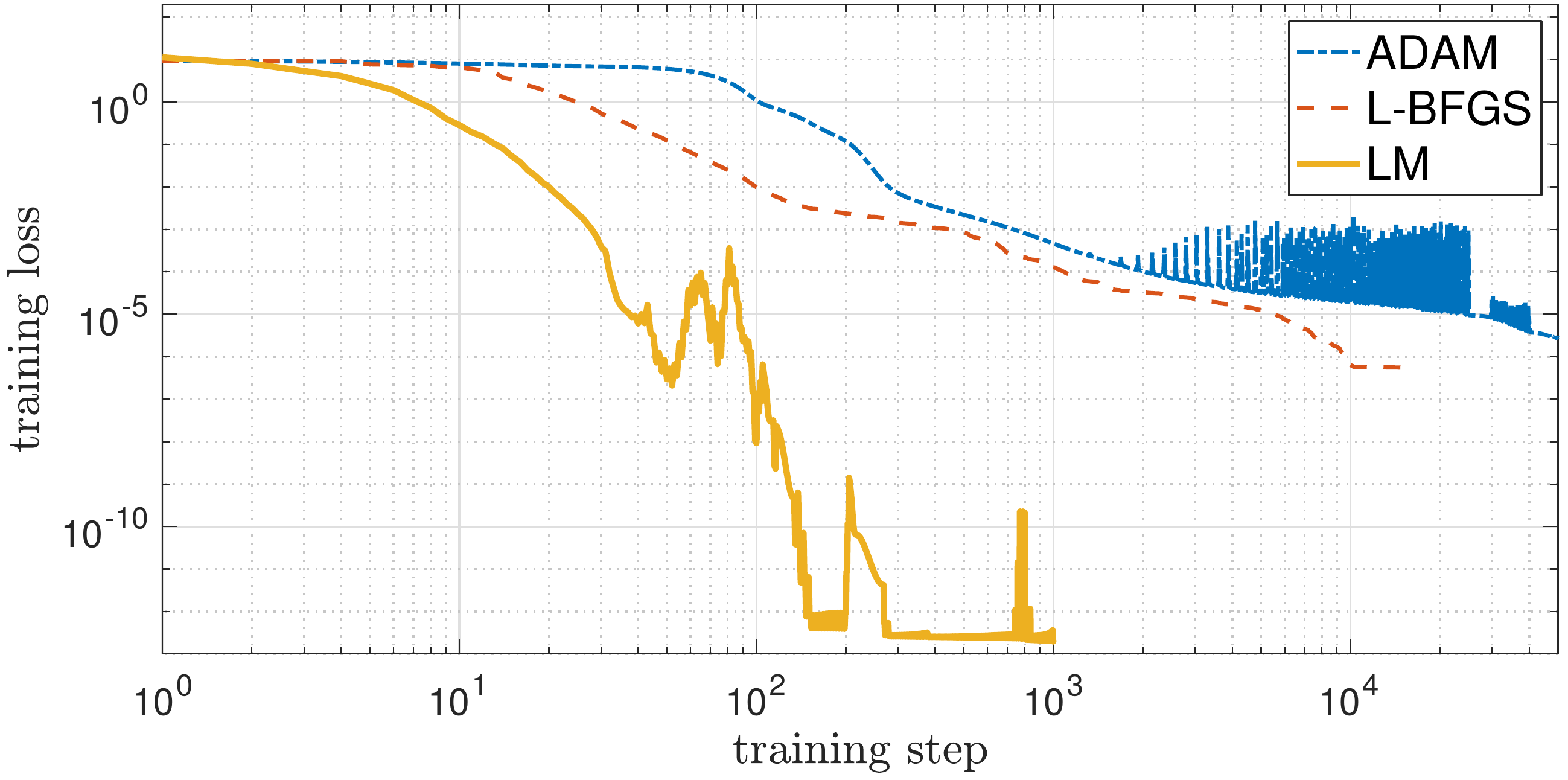}
\caption{Time history of training $\mbox{Loss}_{\Delta_s}$ with $N = 40$ using different optimizers. Dash-dotted line: ADAM, dashed line: L-BFGS, solid line: LM. All training processes use $M = 400$ training points.}
\label{Fig:loss}
\end{figure}

\paragraph{\textbf{Comparisons of loss functions and optimizers.}}
First, we perform the accuracy comparison between our proposed model and existing method in~\cite{FZY21} (i.e., the usage of loss function (\ref{Eq:Loss2})). We fix $M = 400$ training points and train the model using several popular optimizers, such as ADAM~\cite{KB17}, L-BFGS~\cite{LN89}, and LM method.
The results are reported in Table~\ref{Table:loss}, in which the relative $L^2$  errors are shown for $N = 20, 30, 40$ neurons used in the hidden layer.
From the left panel, one can see that the testing accuracy of the present loss model ($\mbox{Loss}_{\Delta_s}$) is quite satisfactory (at least $0.01\%$ predictive accuracy) among all optimizers, showing good approximation capability to the solution for the network model. One can also see that only the results obtained by LM algorithm show convergence tendency with increasing $N$; this is because the LM algorithm, a quadratic convergence method particularly designed for nonlinear least squares problems, generally seeks a local minimum in a faster decaying rate than the other two methods.
As a result, the local minimum found by LM optimizer in general has smaller training loss, and thus achieves higher prediction accuracy. See the time history of training loss for these three optimizers in Fig.~\ref{Fig:loss}.

We also check the testing accuracy using the loss function $\mbox{Loss}_\Delta$ in (\ref{Eq:Loss2}) proposed in~\cite{FZY21}, and show the results  in the right panel of Table~\ref{Table:loss}. One can immediately see how significantly different those relative $L^2$ errors are compared with the results in left panel ($O(10^{-2})$ versus $O(10^{-7})$ for the case of $N=40$ with LM optimizer). And  all errors obtained by $\mbox{Loss}_\Delta$ are apparently greater than $1\%$ no matter which optimizer is adopted. This result indicates that,  the requirement $\partial_n u_\mathcal{N} = 0$ at points along $\Gamma$ gives rise to a locally normal extension solution (in a small neighborhood) which might be complicated, and thus the network model may require more neurons or deeper network structure to be employed to have an accurate prediction.
We further run a series of tests with various exact solutions following the same setup in Table~\ref{Table:loss}.  It turns out that same tendency is observed for both models (not shown here). When other complex surfaces are considered, our model is still able to achieve good predictive accuracy (see later in this subsection) whereas the loss function seeking normal extension by contrast predicts much less accurate solution (these results are not shown here). Based on this finding, we conclude that, with the full expression of differential operators in the loss function, the embedded solution can be accurately expressed under the present shallow neural network.
%For the case of rapidly changing solution on surfaces, we may need to increase the depth or width of neural network structure to attain a desired accuracy, but this goes beyond the scope of the present work and shall be investigated in our future studies.

\begin{table}[h]
\centering
\begin{tabular}{@{}c@{\hspace{15pt}}c@{\hspace{10pt}}c@{\hspace{10pt}}c@{\hspace{10pt}}c@{\hspace{15pt}}c@{\hspace{10pt}}c@{\hspace{10pt}}c@{\hspace{10pt}}}
\hline
\multirow{2}{*}[-3pt]{$N$} & \multicolumn{3}{c}{$\mbox{Loss}_{\Delta_s}$ (\ref{Eq:Loss1}), present work} & & \multicolumn{3}{c}{$\mbox{Loss}_\Delta$ (\ref{Eq:Loss2}), proposed in~\cite{FZY21}} \\
\cmidrule(){2-4}\cmidrule(){6-8}
    & ADAM & L-BFGS & LM &  & ADAM & L-BFGS & LM\\
\midrule
20   & 2.070E$-$04 & 9.500E$-$05 & 6.841E$-$06 & & 7.148E$-$02 & 9.453E$-$02 & 9.774E$-$02 \\
30   & 2.959E$-$04 & 1.009E$-$04 & 1.837E$-$06 & & 5.199E$-$02 & 5.422E$-$02 & 4.390E$-$02 \\
40   & 1.260E$-$04 & 9.376E$-$05 & 3.780E$-$07 & & 4.300E$-$02 & 4.218E$-$02 & 3.304E$-$02 \\
\hline
\end{tabular}
\caption{The average relative $L^2$ errors with different sizes of neurons $N$ in the hidden layer. For each case the number of training points is fixed by $M = 400$.}
\label{Table:loss}
\end{table}

%\paragraph{\textcolor{black}{\textbf{Performance study}}}

\paragraph{\textcolor{black}{\textbf{Comparison of single and double precision computations}}}
\textcolor{black}{Table~\ref{Table:precision} reports an extensive study on the comparison between the single and double precision computations. We vary the number of neurons in the hidden layer and evaluate the relative $L^2$ error between the exact and predicted solutions, terminal loss values, and CPU time (in seconds). In each run, the network model is trained up to $4000$ steps, while the number of training points is fixed by $M = 400$. The results show that for both floating-point representations, given enough training points, the prediction accuracy increases with the number of neurons used.}

%{\bf I am not sure about the spectral accuracy about this paragraph. I would be rather conserved since we only can control the loss but not the error. So I suggest to delete it to make it simple since we have already enough contents.} More specifically, the error decays exponentially as the number of neurons (and the total number of learnable parameters) is doubled. See, for example, the cases $N=5$, $10$ and $20$. This quick convergence can be attributed to the fact that a shallow neural network can approximate a function with spectral accuracy, namely, the error bound is $O(\rho^{-N^{1/d}})$ in $L^{\infty}$-norm~\cite{Mhaskar96}, where $d$ is the dimension of the function variables and $\rho$ is a constant. Similar error bound in $W^{k,\infty}$-norm is obtained for a $\tanh$ activation function~\cite{RLM21}. Therefore, if the loss function and optimizer are chosen properly, such neural network functions can be found and one can expect the prediction accuracy to converge exponentially.

\textcolor{black}{When $N=20$, the loss value obtained using single precision reaches $O(10^{-9})$, which is the limit of single precision calculation. Further increasing the number $N$  does not reduce the loss so the error is stuck at $O(10^{-5})$. For double precision computation, increasing the network complexity beyond $N = 20$ gives a low convergence rate. This is because we stop the training process at $4000$ iterations, the loss may not reach its theoretical minimum. We observe that when the loss is small, its value decays slowly during training, thus requires much more training steps to make the loss smaller.}
%The proper optimizer in this region requires further investigation, which will be left as one of our future work.}

\begin{table}[h]
\centering\color{black}
\begin{tabular}{@{\hspace{5pt}}c@{\hspace{15pt}}c@{\hspace{10pt}}c@{\hspace{10pt}}c@{\hspace{5pt}}c@{\hspace{10pt}}c@{\hspace{10pt}}c@{\hspace{10pt}}c@{\hspace{5pt}}}
\hline
\multirow{2}{*}[-3pt]{$(N, N_p)$} & \multicolumn{3}{c}{double precision} & & \multicolumn{3}{c}{single precision} \\
\cmidrule(){2-4}\cmidrule(){6-8}
    & Error & Loss & CPU time (s) & & Error & Loss & CPU time (s) \\
\midrule
(5, 25)       & 6.524E$-$04 & 1.799E$-$05  & 27   & & 6.967E$-$04 & 2.309E$-$05 & 26 \\
(10, 50)     & 6.237E$-$05 & 7.285E$-$07  & 31   & & 9.746E$-$05 & 6.320E$-$07 & 27 \\
(20, 100)   & 6.841E$-$07 & 4.875E$-$11  & 36   & & 1.349E$-$05 & 4.899E$-$09 & 30 \\
(40, 200)   & 3.780E$-$07 & 1.429E$-$12  & 43   & & 1.314E$-$05 & 2.300E$-$09 & 38 \\
(80, 400)   & 2.082E$-$07 & 1.384E$-$13  & 67   & & 1.398E$-$05 & 2.093E$-$09 & 49 \\
(160, 800) & 1.267E$-$07 & 1.622E$-$13  & 134 & & 9.047E$-$06 & 2.179E$-$09 & 84 \\
\hline
\end{tabular}
\caption{The average relative $L^2$ errors, loss values, and CPU time (in seconds) with different number of neurons $N$ in the hidden layer. For each case the number of training points is fixed by $M = 400$.}
\label{Table:precision}
\end{table}

\paragraph{\textcolor{black}{\textbf{Effect on the number of training points}}}
\textcolor{black}{Next, we investigate the effect on the number of training points. In Table~\ref{Table:LB_training_point}, we deploy $N = 40$ neurons in the hidden layer, and minimize the loss model with the number of training points ranging from $M = 100$ to $M = 500$ (this can be roughly regarded as increasing the spatial resolution in traditional numerical methods). As can be seen, given a small bunch of training data $M = 100$ only results in the accuracy $O(10^{-4})$ with the loss value $O(10^{-8})$. When the loss model is given by the enough information, namely, sufficient number of training points, the network is capable of reaching higher predictive accuracy $O(10^{-7})$ with the loss $O(10^{-12})$.}

\begin{table}[h]
\centering\color{black}
\begin{tabular}{ccc}
\hline
$M$  & Error & Loss\\
\hline
$100$ & 1.423E$-$04 & 3.427E$-$08 \\
$200$ & 1.958E$-$06 & 9.978E$-$11 \\
$300$ & 5.228E$-$07 & 4.777E$-$12 \\
$400$ & 3.780E$-$07 & 1.429E$-$12 \\
$500$ & 2.218E$-$07 & 1.942E$-$13 \\
\hline
\end{tabular}
\caption{The average relative $L^2$ errors for the shallow neural network using different number of $M$. $N = 40$.}
\label{Table:LB_training_point}
\end{table}

\paragraph{\textcolor{black}{\textbf{Effect on the depth of network architecture}}}
\textcolor{black}{We investigate the performance of multiple-hidden-layer network architectures.
With fixed number of training point $M = 400$, we investigate the prediction results using the two-hidden-layer network, which employs $N$ neurons per hidden layer, written as
\begin{align*}
u_\mathcal{N}(\mathbf{x}) = \sum_{j=1}^N \alpha_j \sigma(W^{[2]}_j \sigma(W^{[1]}\mathbf{x}^T+\mathbf{b}^{[1]})+b_j^{[2]}),
\end{align*}
where the weights $\alpha_j\in\mathbb{R}$, $W^{[1]}\in\mathbb{R}^{N\times 3}$ and $W^{[2]}_j\in\mathbb{R}^{1\times N}$, the biases $\mathbf{b}^{[1]}\in\mathbb{R}^N$ and $b_j^{[2]}\in\mathbb{R}$. The total number of learnable parameters is thus counted as $N_p = N^2 + 6N$. From Table~\ref{Table:LB_deep} we can see that,
for the number of learnable parameters $N_p=160, 520, 1840$, the two-hidden-layer network attains equally good accuracy $O(10^{-7})$ in comparison to the shallow one with $N_p=200$, refer the case $N=40, M=400$ in Table~\ref{Table:LB_training_point}. When $N = 80$ is used in each hidden layer, the prediction accuracy reaches $O(10^{-8})$ and the loss value decays as low as $O(10^{-14})$. Whereas this small improvement of accuracy requires a large number of parameters $N_p = 6880$ needed to be trained. Thus, the usage of shallow neural network representation is readily able to encode smooth solutions, and the accuracy performance is equally well compared to the two-hidden-layer network.}

\begin{table}[h]
\centering\color{black}
\begin{tabular}{ccc}
\hline
$(N, N_p)$  & Error & Loss\\
\hline
$(10, 160)$ & 2.347E$-$07 & 2.054E$-$12 \\
$(20, 520)$ & 3.418E$-$07 & 2.617E$-$13 \\
$(40, 1840)$ & 1.316E$-$07 & 1.170E$-$13 \\
$(80, 6880)$ & 5.445E$-$08 & 2.582E$-$14 \\
\hline
\end{tabular}
\caption{The average relative $L^2$ errors for the two-hidden-layer neural network. In this case the total number of learnable parameters is counted by $N_p = N^2+6N$ and the training data points is fixed by $M = 400$.}
\label{Table:LB_deep}
\end{table}

\textcolor{black}{\paragraph{\textbf{Application to a non-closed surface}}
When the considered surface is not closed, the underlying PDE must be subject to an additional boundary condition along $\partial\Gamma$.
Here, we consider the Dirichlet-type boundary condition $u(\bx) = u_b(\bx)$ for $\bx\in\partial\Gamma$, so,
it is straightforward to simultaneously enforce mean squared errors for both differential equation and boundary condition in a loss function. That is, given training sets $\{\mathbf{x}^i\in\Gamma\}_{i=1}^M$ and $\{\mathbf{x}^j_{\partial\Gamma}\in\partial\Gamma\}_{j=1}^{M_b}$, the loss function~(\ref{Eq:Loss1}) is thus slightly modified with an additional penalty term as
\begin{align*}
\mbox{Loss}_{\Delta_s}(\mathbf{p}) &= \frac{1}{M}\sum_{i=1}^M \left[\Delta u_\mathcal{N}(\mathbf{x}^i) - 2H(\mathbf{x}^i)\partial_nu_\mathcal{N}(\mathbf{x}^i) - \mathbf{n}(\mathbf{x}^i)^T\left(\nabla^2u_\mathcal{N}(\mathbf{x}^i)\right)\mathbf{n}(\mathbf{x}^i) - f(\mathbf{x}^i)\right]^2 \\
&+ \frac{1}{M_b}\sum_{j=1}^{M_b}\left[ u_\mathcal{N}(\mathbf{x}^j_{\partial\Gamma}) - u_b(\mathbf{x}^j_{\partial\Gamma}) \right]^2.
\end{align*}
We run a test example whose solution is chosen as $u(x,y,z) = \sin(x)\exp(\cos(y-z))$ and the hemi-elliposid $\psi(x,y,z) = (x/1.5)^2 + (y)^2 + (z/0.5)^2 - 1$ with $z>0$ is considered.
%Since the above loss model takes the form of least squares error, it can be trained efficiently using the LM optimizer.
In Table~\ref{Table:LB_boundary}, with fixed number of training points $M = 400$ and $M_b = 100$, we investigate the prediction accuracy with different number of neurons $N$ in the hidden layer.  As seen, with the presence of boundary conditions, the proposed model is still able to attain satisfactory accuracy. Again, given a sufficient number of training points $M$ and $M_b$, the prediction accuracy increases as the number $N$  increases.}

\begin{table}[h]
\centering\color{black}
\begin{tabular}{cc}
\hline
$(N,N_p)$    & Error \\
\hline
$(5,25)$       & 8.903E$-$04  \\
$(10,50)$     & 6.598E$-$05  \\
$(20,100)$   & 1.828E$-$06  \\
$(40,200)$   & 1.735E$-$07  \\
$(80,400)$   & 8.158E$-$08  \\
$(160,800)$ & 7.697E$-$08  \\
\hline
\end{tabular}
\caption{The average relative $L^2$ errors. For each case the number of training points is fixed by $M = 400$ and $M_b = 100$.}
\label{Table:LB_boundary}
\end{table}

\paragraph{\textbf{Numerical results of more complex surfaces}}
In the previous tests,  we only focus on the surface geometry as simple as an ellipsoid (or hemi-ellipsoid).  Here, we present the numerical accuracy results for our proposed neural network using the loss function (\ref{Eq:Loss1}) with more complex geometries such as torus, genus-2 surface, and cheese-like surface (see Fig.~\ref{Fig:curvature}). %The LM method is used as the optimizer. %the training procedure is finished when the mean squared loss is less than $10^{-11}$ or 1000 training steps is reached.

In Table~\ref{Table:LB} we show the average $L^2$ relative errors for those different surfaces.  Again,  we choose $u(x,y,z) = \sin(x)\exp(\cos(y-z))$ and fix $M = 400$ training points which are randomly deployed along each surface, and vary the number of neurons used in the hidden layer from $N = 20,30,40,50,60$. 
We see that for all those different surface geometries, using just $N = 20$ neurons (learnable parameters $N_p = 100$) is sufficient to encode the solutions with at least $0.01\%$ predictive accuracy. Although the numerical convergence is not rigorously verified, the increase of neurons generally leads to better accuracy for all these cases shown in the table.

\begin{table}[h]
\centering
\begin{tabular}{c|ccc}
\hline
$(N,N_p)$ & torus & genus-2 & cheese-like \\
\hline
$(20,100)$ & 2.774E$-$05  & 1.816E$-$06 & 1.522E$-$04 \\
$(30,150)$ & 5.568E$-$06  & 9.150E$-$07 & 2.897E$-$05 \\
$(40,200)$ & 2.181E$-$06  & 6.100E$-$07 & 1.018E$-$05 \\
$(50,250)$ & 1.708E$-$06  & 4.731E$-$07 & 7.176E$-$06 \\
$(60,300)$ & 1.139E$-$06  & 5.169E$-$07 & 5.617E$-$06 \\
\hline
\end{tabular}
\caption{The average relative $L^2$ errors with different number of neurons $N$ in the  hidden layer. For each case the number of training points is fixed by $M = 400$.}
\label{Table:LB}
\end{table}

\begin{figure}[!]
\centering
\includegraphics[scale=0.41]{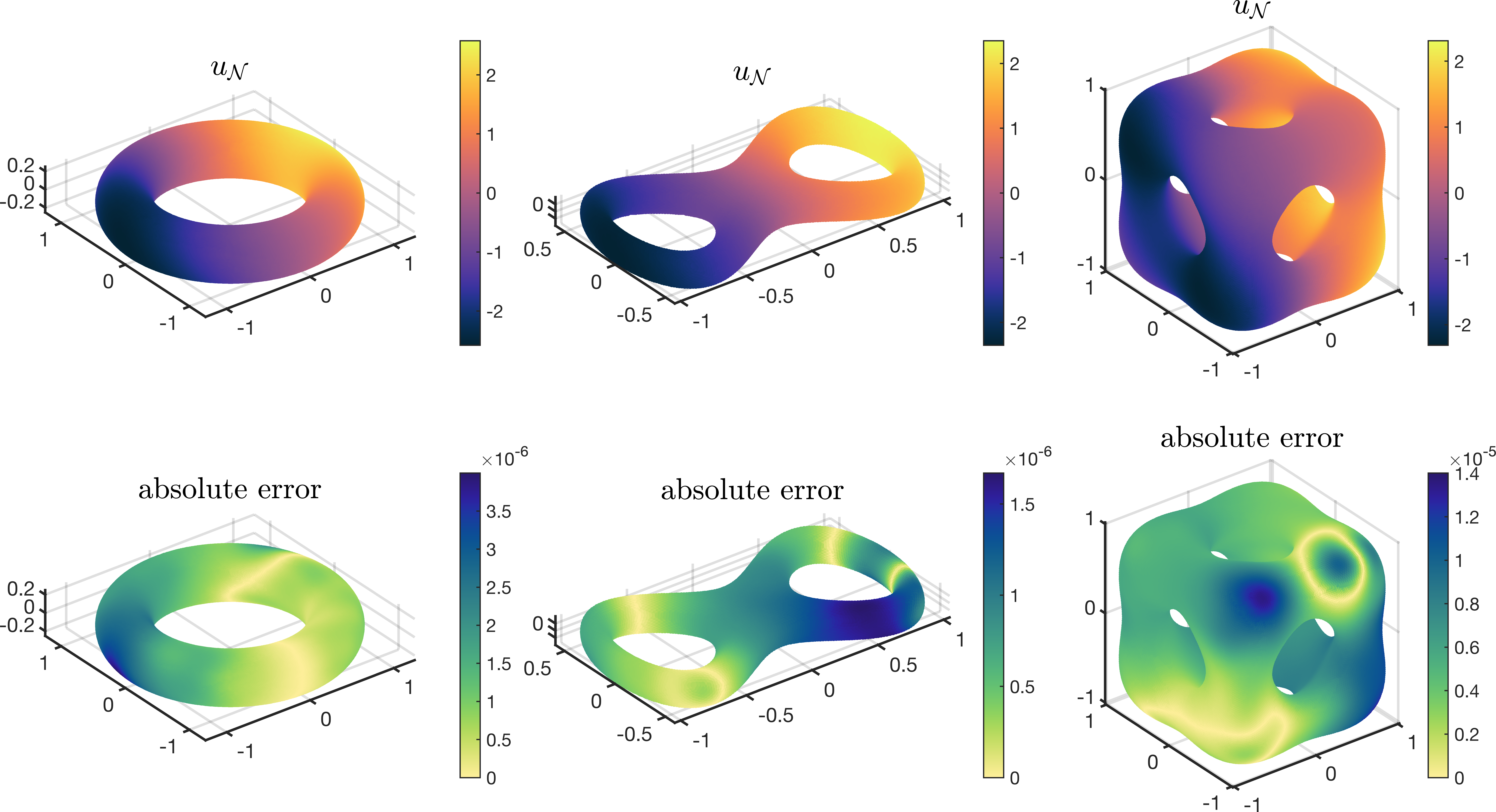}
\caption{Prediction solution $u_\mathcal{N}$ and corresponding absolute error $|u_\mathcal{N}-u|$ with $N = 60$ neurons employed. From left to right: torus, genus-2 surface, cheese-like surface.}
\label{Fig:error_LB}
\end{figure}

The predicted network solution $u_\mathcal{N}$ (with $N = 60$ and $M = 400$) and the absolute error $|u_\mathcal{N}-u|$ for these surfaces are depicted in Fig.~\ref{Fig:error_LB}. One can see that, regardless of the surface geometries, our designed network model is able to obtain equally accurate prediction for all cases (the largest absolute error does not necessarily occur at high curvature points). In addition, these results are obtained by randomly sampled training points on the underlying surfaces, highlighting the robustness feature of the mesh-free nature of the neural network model.

%%%%%%%%%%%%%%%%%%%%%%%%%%%%%%%%%%%%%%%%%
\section{Time-dependent PDEs on static surfaces}\label{Sec:time_PDE}
%%%%%%%%%%%%%%%%%%%%%%%%%%%%%%%%%%%%%%%%%

In this section, we turn our attention to solve time-dependent PDEs on static surfaces. Given a regular and closed surface $\Gamma$, along which we consider the PDEs of the general form
\begin{align}\label{Eq:PDE2}
\partial_t u(\mathbf{x},t) = \mathcal{L}(u(\mathbf{x},t)) + f(\mathbf{x},t) \quad \mbox{ on } \Gamma,\, t\in(0,T],
\end{align}
where $t$ denotes the time variable and $T$ is the terminal time; $f$ is a source term defined on $\Gamma$. Again $\mathcal{L}(u)$ may contain the surface gradient $\nabla_s u$, surface diffusion $\Delta_s u$, or $\nabla_s\cdot\mathbf{v}$ for some known vector field $\mathbf{v}$.
%Although $\mathcal{L}$ is not complete, it consists of most surface differential terms of practical interests in applications.

In this section, we shall concentrate on solving the surface diffusion equation ($\mathcal{L} = \Delta_s$) as
\begin{align}\label{Eq:diffusion}
\partial_t u(\mathbf{x},t) = \Delta_s u(\mathbf{x},t) + f(\mathbf{x},t) \quad \mbox{ on } \Gamma,\, t\in(0,T].
\end{align}
The above PDE is subjected to an initial condition
\begin{align}\label{Eq:IC}
u(\mathbf{x},t=0) = u_0(\mathbf{x}) \quad \mbox{ on } \Gamma.
\end{align}
To solve this time-dependent PDE, we follow the pioneer framework of physics-informed neural networks proposed in~\cite{raissi2019}, i.e., the above diffusion equation is solved by continuous-time or discrete-time model.
%\revone{In addition, we propose an initial free learning method in which the obtained solution fulfills the initial condition exactly. These methods are stated as followings.}

\paragraph{\textbf{Continuous-time model}}
It is natural to encapsulate both spatial and time variables as the input of neural network function.
Thus, the approximate solution to Eq.~(\ref{Eq:diffusion}) now can be written as
\begin{align}\label{NN:cont_time}
u_\mathcal{N}(\mathbf{x},t) = \sum_{j=1}^N \alpha_j \sigma(W_j(\mathbf{x},t)^T+b_j).
\end{align}
Differing from the stationary case (see Eq.~(\ref{Eq:NN})) due to the time variable augmentation, the dimension of weights becomes $W_j\in\mathbb{R}^{1\times4}$, so the total number of learnable parameters is increased as $N_p = 6N$.

To learn those parameters, as in stationary case, we train the neural net model using the identities in Eq.~(\ref{Eq:diff_operator}) to compute differential terms appeared in Eq.~(\ref{Eq:diffusion}). Thus,
it is straightforward to employ the physics-informed learning method to minimize the mean squared residual for both differential equation (\ref{Eq:diffusion}) and initial condition (\ref{Eq:IC}). For given training points $\{(\mathbf{x}^i,t^i)|\mathbf{x}^i\in\Gamma, t^i\in(0,T]\}_{i=1}^{M_T}$ and $\{\mathbf{x}_0^j\in\Gamma\}_{j=1}^{M_0}$, the natural training loss is chosen as
\begin{align}\label{Eq:Loss_cont_time}
\mbox{Loss}_c(\mathbf{p}) = \frac{1}{M_T}\sum_{i=1}^{M_T} \left[\partial_t u_\mathcal{N}(\mathbf{x}^i,t^i) - \Delta_s u_\mathcal{N}(\mathbf{x}^i,t^i) - f(\mathbf{x}^i,t^i)\right]^2 + \frac{1}{M_0}\sum_{j=1}^{M_0}\left[ u_\mathcal{N}(\mathbf{x}^j_0,0) - u_0(\mathbf{x}^j_0) \right]^2,
\end{align}
where $\Delta_s u_\mathcal{N}(\mathbf{x}^i,t^i) = \Delta u_\mathcal{N}(\mathbf{x}^i,t^i) - 2H(\mathbf{x}^i)\partial_nu_\mathcal{N}(\mathbf{x}^i,t^i) - \mathbf{n}(\mathbf{x}^i)^T\left(\nabla^2u_\mathcal{N}(\mathbf{x}^i,t^i)\right)\mathbf{n}(\mathbf{x}^i)$.

\paragraph{\textbf{Discrete-time model}}
In contrast to the continuous-time model, in discrete-time model the PDE (\ref{Eq:diffusion}) is alternatively solved by a semi-discretization scheme as in classical numerical methods~\cite{raissi2019}. That is, we obtain the numerical solution via the $q$-stage time-stepping implicit Runge-Kutta (RK) scheme:
\begin{align}
& u^{n+c_j} = u^n + \Delta t \sum_{k=1}^q a_{jk}( \Delta_s u^{n+c_k} + f^{n+c_k} ), \quad j = 1,2,\cdots,q, \label{Eq:RK}\\
& u^{n+1} = u^n + \Delta t \sum_{k=1}^q b_k ( \Delta_s u^{n+c_k} + f^{n+c_k} ), \label{Eq:update}
\end{align}
where $\Delta t$ is the time step size, $u^{n+c_j} = u(\mathbf{x}, (n+ c_j)\Delta t)$ and $f^{n+c_j} = f(\mathbf{x}, (n+ c_j)\Delta t)$ are the intermediate solution and source term correspondingly, and $u^{n+1} = u(\mathbf{x},(n+1)\Delta t)$ is the numerical solution at the next time level. Here we adopt Gauss-Legendre method so the temporal discretization error of above $q$-stage Runge-Kutta scheme is $O(\Delta t^{2q})$, where the parameters $\{a_{jk},b_k,c_k\}$ are given from Butcher tableau~\cite{Iserles09}. By taking sufficiently large $q$, this high-order scheme allows us to obtain an accurate numerical solution $u^{n+1}$ even with large $\Delta t$.  Meanwhile, the numerical stability can be retained due to the full implicity in Eq.~(\ref{Eq:RK}).

To obtain $u^{n+1}$, we need to learn those intermediate network solutions, $u_\mathcal{N}^{n+c_j}$, again via physics-informed learning technique.
We proceed by placing a multi-output neural network $\mathbf{u}_\mathcal{N}(\mathbf{x}) = [u_\mathcal{N}^{n+c_1}(\mathbf{x}), u_\mathcal{N}^{n+c_2}(\mathbf{x}), \cdots, u_\mathcal{N}^{n+c_q}(\mathbf{x}), u_\mathcal{N}^{n+1}(\mathbf{x})]^T$ and it can be compactly expressed by
\begin{align*}
\mathbf{u}_\mathcal{N}(\mathbf{x}) = W^{[2]}\sigma(W^{[1]}\mathbf{x}^T+\mathbf{b}^{[1]}),
\end{align*}
where $W^{[1]}\in\mathbb{R}^{N\times3}$ and $W^{[2]}\in\mathbb{R}^{(q+1)\times N}$ are the weight matrices and $\mathbf{b}^{[1]}\in\mathbb{R}^N$ is the bias (so all $u_\mathcal{N}^{n+c_j}$ and $u_\mathcal{N}^{n+1}$ are learned in a single network).
In this network, there are $N_p = (5+q)N$ parameters needed to be learned. The loss function is thereby designed to simultaneously enforce all discretization equations (\ref{Eq:RK}) together with the updating step (\ref{Eq:update}). That is,  given a set of training points $\{\mathbf{x}^i\in\Gamma\}_{i=1}^M$, we have
\begin{align}\label{Eq:Loss_discrete_time}
\begin{split}
\mbox{Loss}_d(\mathbf{p}) &= \frac{1}{M}\sum_{j=1}^q\sum_{i=1}^M \left[ u_\mathcal{N}^{n+c_j}(\mathbf{x}^i) - u^n(\mathbf{x}^i) - \Delta t \sum_{k=1}^q a_{jk}( \Delta_s u_\mathcal{N}^{n+c_k}(\mathbf{x}^i) + f^{n+c_k}(\mathbf{x}^i) ) \right]^2\\
&+ \frac{1}{M}\sum_{i=1}^M \left[ u_\mathcal{N}^{n+1}(\mathbf{x}^i) - u^n(\mathbf{x}^i) - \Delta t \sum_{k=1}^q b_{k}( \Delta_s u_\mathcal{N}^{n+c_k}(\mathbf{x}^i) + f^{n+c_k}(\mathbf{x}^i) ) \right]^2.
\end{split}
\end{align}
After finishing the training of the above model loss, we then use this prediction $u_\mathcal{N}^{n+1}$ as the initial condition to advance to the next time level $u^{n+2}_\mathcal{N}$ by proceeding to the same training process.  Eventually, we obtain the numerical solution at the target terminal time.
%We should point out that, in \cite{raissi2019} the numerical solution $u^{n+1}$ is also included in the multi-output network and thus both Eq(\ref{Eq:RK})-(\ref{Eq:update}) are penalized in a training loss, whereas we note that the updating step (\ref{Eq:update}) can be done explicitly with the given intermediate solutions. Comparing with their original model, the present model has less parameters to be trained and thus it indeed saves some computational efforts.

%%%%%%%%%%%%%%%%%%%%%%%%%%%%%%%%%%%%%%%%%
\subsection{Numerical accuracy}
%%%%%%%%%%%%%%%%%%%%%%%%%%%%%%%%%%%%%%%%%

We perform the capability of continuous- and discrete-time neural network model, corresponding to $\mbox{Loss}_c$ in (\ref{Eq:Loss_cont_time}) and $\mbox{Loss}_d$ in (\ref{Eq:Loss_discrete_time}), for encoding the diffusion equation on the cheese-like surface.
We check the prediction accuracy by considering the exact solution
\begin{align*}
u(x,y,z,t) = \sin(x+\sin(t))\exp(\cos(y-z)),
\end{align*}
so the source term $f$ can be obtained accordingly. We set the  terminal time $T = 1$.
For continuous-time model we use $M_0 = 100$ and $M_T = 800$ spatial-temporal training points, in which the surface points $\mathbf{x}^i$ are randomly sampled while temporal points $t^i$ are chosen based on Latin Hypercube Sampling strategy~\cite{Stein87}.
%The same number $M_T = 800$ is also used in the initial free model.
In discrete-time model we set $M = 200$ spatial training points and adopt 6-stage implicit Runge-Kutta scheme with time step size $\Delta t = 1$ (so the network solution at terminal time $T = 1$ is obtained under a single time update).
The average relative $L^2$ errors at $T = 1$ for network models with various neurons of the hidden layer $N$ are shown in Table~\ref{Table:time}.
Again, both models can obtain accurate predictive results.
%As one can see, the prediction accuracy of the initial free loss performs better than the one obtained by the usual continuous-time model (especially for small $N$).
Furthermore, as expected, the increase of the number of neurons generally leads to better accuracy as well.

\begin{table}[h]
\centering
\begin{tabular}{cc|cc}
\hline
$(N,N_p)$ & continuous-time model  & $(N,N_p)$ & discrete-time model \\
\hline
$(20,120)$ & 1.400E$-$03 & $(20,220)$ & 6.448E$-$04 \\
$(30,180)$ & 1.975E$-$04  & $(30,330)$ & 5.013E$-$05 \\
$(40,240)$ & 1.390E$-$04  & $(40,440)$ & 1.627E$-$05 \\
$(50,300)$ & 5.984E$-$05  & $(50,550)$ & 7.920E$-$06 \\
$(60,360)$ & 3.661E$-$05  & $(60,660)$ & 6.446E$-$06 \\
\hline
\end{tabular}
\caption{The average relative $L^2$ errors at $T = 1$ for continuous-time and discrete-time model with various neurons of hidden layer $N$.
In each test, we fix $M_T = 800$ and $M_0 = 100$ for the continuous-time model; $M = 200$ and $\Delta t = 1$ for the 6-stage RK discrete-time model.}
\label{Table:time}
\end{table}

%%%%%%%%%%%%%%%%%%%%%%%%%%%%%%%%%%%%%%%%%
\subsection{A surface heating up application}
%%%%%%%%%%%%%%%%%%%%%%%%%%%%%%%%%%%%%%%%%

\begin{figure}[h]
\centering
\includegraphics[scale=0.45]{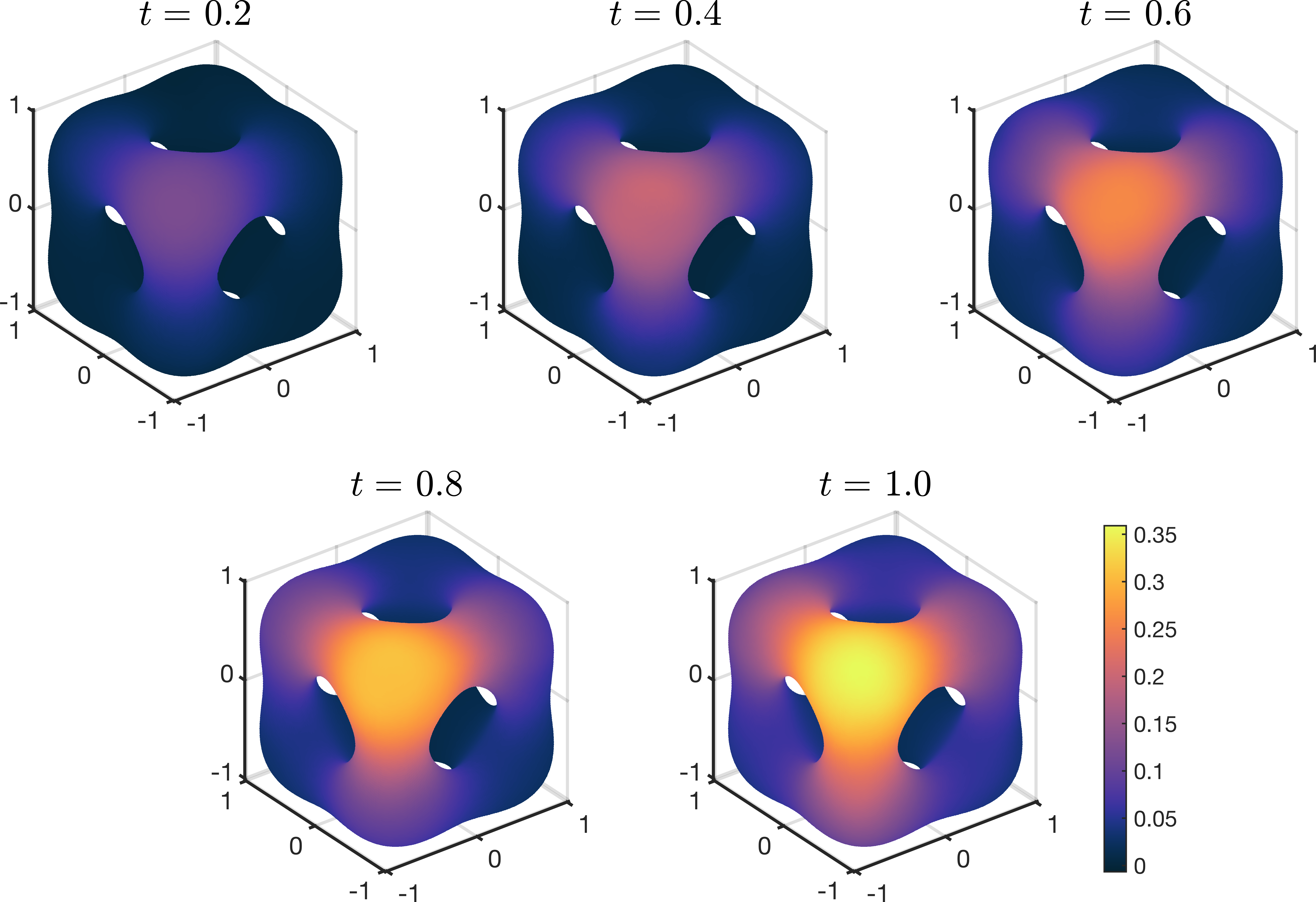}
\caption{The snapshots of solution distribution for heating up the cheese-like surface at different times. The color code ranging from 0 to 0.36 indicates the magnitude of the solution.}
\label{Fig:heating}
\end{figure}

We perform an application simulation by mimicking the process of heating up a surface. The initial condition is set to be zero everywhere and the heating source is a time-independent Gaussian bump given by
\begin{align*}
f(x,y,z) = \exp(-( (x+1)^2 + (y+1)^2 + (z-1)^2 )),
\end{align*}
so the majority of the source accumulates in the vicinity of the point $(-1,-1,1)$. The diffusion equation is solved using the discrete-time model with the 4-stage RK scheme, in which we set the time step $\Delta t = 0.1$ and compute the solution up to the terminal time $T = 1$. We use $M = 500$ training points and $N = 100$ neurons. \textcolor{black}{Since there is no analytical solution available in this case, we are unable to measure the relative $L^2$ error quantitatively. We simply train the network to get  the loss value  to the order of magnitude $O(10^{-8})$ which is roughly matched with the temporal discretization error $(\Delta t)^8 = 10^{-8}$.
The successive snapshots of time-evolutional solution are displayed in Fig.~\ref{Fig:heating}. As we can see, near the source of Gaussian bump, the magnitude of prediction solution becomes larger as time evolves. At the same time, the  heat distribution becomes wider  due to the diffusion mechanism in the PDE model. Therefore, the predictive solution generated by our network model presents some visually plausible results.}

%%%%%%%%%%%%%%%%%%%%%%%%%%%%%%%%%%%%%%%%%
\section{PDEs on evolving surfaces}\label{Sec:evolving_PDE}
%%%%%%%%%%%%%%%%%%%%%%%%%%%%%%%%%%%%%%%%%

In this section, we extend the proposed methodology to solve PDEs on evolving surfaces. Here, the considered surface $\Gamma(t)$ evolves with a prescribed velocity field $\mathbf{v}(\bx,t)$ so %$\Gamma = \Gamma(t) = \{\bx\in\mathbb{R}^3\}$, and
its configuration follows the evolutional equation
\begin{align}\label{Eq:surface_evolve}
\partial_t\bx = \mathbf{v}(\bx(t),t), \quad \bx(t) \in \Gamma(t),\, t\in(0,T]
\end{align}
together with an initial configuration $\Gamma(0)$ represented by $\bx(t=0) = \bx_0$. For simplicity, we assume that $\Gamma(t)$ remains a regular surface under the velocity field $\mathbf{v}$.

%Due to the presence of the transportation velocity, the time-dependent PDE on $\Gamma(t)$, in Eulerian coordinates (thus $u = u(\bx,t)$), takes the general from
%\begin{align*}
%\partial_tu + \mathbf{v}\cdot\grad u + (\grad_s\cdot\mathbf{v})u = \mathcal{L}(u) + f \quad \mbox{ on } \Gamma(t),\, t\in(0,T],
%\end{align*}
%where $\partial_tu + \mathbf{v}\cdot\grad u$ and $(\grad_s\cdot\mathbf{v})u$ denote the material derivative and surface stretching term respectively, and $f = f(\bx,t)$ is a given source term defined on $\Gamma(t)$.
Throughout this section, we consider the following advection-diffusion equation on the evolving surface $\Gamma(t)$ as
\begin{align}\label{Eq:advection_diffusion}
\partial_tu + \mathbf{v}\cdot\grad u + (\grad_s\cdot\mathbf{v})u = \Delta_s u + f \quad \mbox{ on } \Gamma(t),\, t\in(0,T],
\end{align}
where $\partial_tu + \mathbf{v}\cdot\grad u$ denotes the material derivative of $u$,  and the term $(\grad_s\cdot\mathbf{v})u$ represents the surface stretching effect on the quantity $u$. The term $f = f(\bx,t)$ is again a given source term defined on $\Gamma(t)$. Of course, the above
equation must be accompanied with a given initial condition $u(\bx_0,t = 0) = u_0(\bx_0)$ on $\Gamma(t=0)$. One should note that the above equation
(\ref{Eq:advection_diffusion}) is popularly used in modeling certain physical applications; for instance,  the insoluble surfactant concentration along a droplet surface in fluid flows \cite{SHL18, HCLT19}. The major challenge of solving Eq.~(\ref{Eq:advection_diffusion}) arises from the time-dependent computation of surface geometrical quantities such as mean curvatures $H(\bx,t)$ and normal vectors along the evolving surface (which are involved in those conventional differential terms as seen in Eq.~(\ref{Eq:diff_operator})). We aim to solve the PDE system~(\ref{Eq:surface_evolve})-(\ref{Eq:advection_diffusion}) under a unified continuous-time neural network framework, as stated as follows.

%%%%%%%%%%%%%%%%%%%%%%%%%%%%%%%%%%%%%%%%%
\subsection{Neural network solver for PDEs on evolving surfaces}
%%%%%%%%%%%%%%%%%%%%%%%%%%%%%%%%%%%%%%%%%

To track this time evolving surface using neural network representation, we adopt the surface parametrization as
\begin{align*}
\Gamma(t) = \{\bx(\theta, \phi, t)\in\mR^3\,| \, (\theta, \phi)\in[0, \pi]\times[0, 2\pi), t\in[0,T]\}.
\end{align*}
The key observation is that a closed surface (with genus zero) is homeomorphic to a unit sphere $\mS^2$ in three-dimensional space. Therefore, there exists a continuous and invertible mapping between the surface $\Gamma(t)$ and $\mS^2$. We hereby propose a two-hidden-layer neural network structure to represent the surface. We first map the input variables $(\theta, \phi)$ to the unit sphere $\mS^2$ (as the output of the first hidden layer), and then use a fully-connected neural network to learn the homeomorphism between $\mS^2$ and the surface. More precisely, let $S^2(\theta, \phi) = (\sin\theta\cos\phi, \sin\theta\sin\phi, \cos\theta)$, then the homeomorphic network can be written as
\begin{align}\label{Eq:NN_surface}
\bx_\mathcal{N}(\theta, \phi, t) = W^{[2]}\sigma(W^{[1]}(S^2(\theta, \phi), t)^T + \mathbf{b}^{[1]}),
\end{align}
where the weight matrices $W^{[1]}\in\mathbb{R}^{N\times4}$ and $W^{[2]}\in\mathbb{R}^{3\times N}$, and the bias $\mathbf{b}^{[1]}\in\mathbb{R}^N$ (so the total number of training parameters $N_p = 8N$).
We should emphasize the features of the above surface network representation:
(i) There are no learnable parameters needed to be trained from the input layer to the first hidden layer, i.e., the output of the first hidden layer is directly computed through the nonlinear map $S^2$, whereas, the remaining parameters of the network ($W^{[1]}$, $W^{[2]}$, and $\mathbf{b}^{[1]}$) need to be trained.
(ii) This representation automatically fulfills $2\pi$-periodicity in $\phi$-direction while the pole conditions at $\theta = 0$ and $\pi$ are taken care by the parametrization of $S^2$.
(iii) When a genus $g$ surface is considered, following the same idea, we can adopt the mapping from the parametric domain to a $g$-torus so that the homeomorphism can be learned using neural network representation.

Now, with the proper surface network representation~(\ref{Eq:NN_surface}), we proceed to solve the surface evolving equation~(\ref{Eq:surface_evolve}). Given sets of training points $\{(\theta^i,\phi^i,t^i) | (\theta^i, \phi^i)\in[0, \pi]\times[0, 2\pi), t^i\in(0,T]\}_{i=1}^{M_T}$ and $\{(\theta_0^j, \phi_0^j)\in[0, \pi]\times[0, 2\pi)\}_{j=1}^{M_0}$, the surface configuration at any instantaneous time is found by minimizing the continuous-time loss model as
\begin{align}\label{Eq:Loss_surface}
\mbox{Loss}_\bx(\bp) = \frac{1}{M_T}\sum_{i=1}^{M_T}\left[ \partial_t\bx_\mathcal{N}(\theta^i,\phi^i,t^i) - \mathbf{v}(\bx_\mathcal{N}(\theta^i,\phi^i,t^i),t^i) \right]^2
+\frac{1}{M_0}\sum_{j=1}^{M_0}\left[ \bx_\mathcal{N}(\theta_0^j,\phi_0^j,0) - \bx_0(\theta_0^j,\phi_0^j) \right]^2.
\end{align}
Here, both $(\theta^i,\phi^i)$ and $(\theta_0^j,\phi_0^j)$ are chosen so that $S^2$ acting on those points are randomly distributed on $\mS^2$.
This strategy shall effectively avoid local cluster of sample points on $\Gamma(t)$. After the termination of the training process, we use the network solution $\bx_\mathcal{N}(\theta,\phi,t)$ to build up the training sets $\{(\bx^i,t^i) | \bx^i\equiv \bx_\mathcal{N}(\theta^i,\phi^i,t^i)\in\Gamma(t^i)\}_{i=1}^{M_T}$ and $\{\bx_0^j\equiv\bx_\mathcal{N}(\theta_0^j,\phi_0^j,0)\in\Gamma(0)\}_{j=1}^{M_0}$, and find the normal vectors $\mathbf{n}(\bx^i,t^i)$ and mean curvatures $H(\bx^i,t^i)$ via the first and second fundamental forms in differential geometry~\cite{Walker15}. As a consequence, finding the solution to the surface PDE~(\ref{Eq:advection_diffusion}) is a straightforward application of the present method. Namely, expressing the shallow neural network solution by Eq.~(\ref{NN:cont_time}), we minimize the loss function
\begin{align}\label{Eq:Loss_PDE_evolve}
\begin{split}
\mbox{Loss}_u(\mathbf{p}) &= \frac{1}{M_T}\sum_{i=1}^{M_T} \left[\partial_t u_\mathcal{N}(\mathbf{x}^i,t^i) + \mathbf{v}(\bx^i,t^i)\cdot\grad u(\bx^i,t^i) + (\grad_s\cdot\mathbf{v}(\bx^i,t^i))u(\bx^i,t^i) - \Delta_s u_\mathcal{N}(\mathbf{x}^i,t^i) - f(\mathbf{x}^i,t^i)\right]^2 \\&+ \frac{1}{M_0}\sum_{j=1}^{M_0}\left[ u_\mathcal{N}(\mathbf{x}^j_0,0) - u_0(\mathbf{x}^j_0) \right]^2,
\end{split}
\end{align}
where $\Delta_s u_\mathcal{N}(\mathbf{x}^i,t^i) = \Delta u_\mathcal{N}(\mathbf{x}^i,t^i) - 2H(\mathbf{x}^i,t^i)\partial_nu_\mathcal{N}(\mathbf{x}^i,t^i) - \mathbf{n}(\mathbf{x}^i,t^i)^T\left(\nabla^2u_\mathcal{N}(\mathbf{x}^i,t^i)\right)\mathbf{n}(\mathbf{x}^i,t^i)$, and $\nabla_s\cdot\mathbf{v}(\bx^i,t^i) = \nabla\cdot\mathbf{v}(\bx^i,t^i) - \mathbf{n}(\bx^i,t^i)^T\grad\mathbf{v}(\bx^i,t^i)\mathbf{n}(\bx^i,t^i)$.

It is worth mentioning that, in other Eulerian coordinates based embedding methods~\cite{LLZ11, HCLT19}, they require an operator splitting strategy so that the advection and diffusion parts are solved separately in order to find the solution. In comparison, the present method~(\ref{Eq:Loss_PDE_evolve}) deals with the surface PDE at the instantaneous time $t = t_i$ directly; thus, the implementation is simple and straightforward.

Since the surface configuration $\bx_\mathcal{N}$ and the underlying solution $u_\mathcal{N}$ change simultaneously as
time proceeds, it is more practical to obtain them in a time sequential manner especially for longer time $T$. In this way, we divide the time interval $[0, T]$ into $n$ uniform subintervals as $[0, T]=\cup_{k=1}^{n} [T_{k-1}, T_k]$, and apply the above learning machinery to obtain the solutions of $\bx_\mathcal{N}$ and $u_\mathcal{N}$ in each time interval $[T_{k-1}, T_k]$ starting at $k=1$. We repeatedly use the loss functions
Eq.~(\ref{Eq:Loss_surface}) and Eq.~(\ref{Eq:Loss_PDE_evolve}) by resuming the initial data that is obtained from the trained results in the previous time interval. Unless otherwise stated, we use the notation $M$ (instead of $M_T$) to denote the number of training points used in parametric domain $(\theta, \phi,t)\in[0, \pi]\times[0, 2\pi)\times (T_{k-1},T_k]$, and $M_0$ to denote the number of training points used in the initial data at $T_{k-1}$.  Here, we randomly choose $M_{test} = 100M$ test points and repeat the numerical experiments 5 times
so the average relative $L^2$ error is computed based on these 5 runs on different training points. The numerical results are shown in the following subsections.

%%%%%%%%%%%%%%%%%%%%%%%%%%%%%%%%%%%%%%%%%
\subsection{Numerical results}
%%%%%%%%%%%%%%%%%%%%%%%%%%%%%%%%%%%%%%%%%
This example aims to demonstrate the capability and accuracy of the proposed method for solving the advection-diffusion equation on 2D evolving surface in $\mR^3$.
We consider the case of an oscillating ellipsoid~\cite{PR16, HCLT19} whose configuration is described by
\begin{align*}
\Gamma(t) = \left\{(x,y,z)\Bigg| \left(\frac{x}{1.5a(t)}\right)^2 + y^2 + \left(\frac{z}{0.5}\right)^2 = 1 \right\}.
\end{align*}
The associated velocity field is $\mathbf{v} = \left(\frac{a'(t)}{a(t)}x, 0, 0\right)$ and we set $a(t) = \sqrt{1+0.95\sin(\pi t)}$.
One should note that the above Cartesian representation for $\Gamma(t)$ can be easily rewritten as the parametric form in terms of $(\theta, \phi)$
using unit sphere representation $S^2$ in previous subsection.
The exact solution to the surface advection-diffusion equation (\ref{Eq:advection_diffusion}) is again chosen as $u(x,y,z,t) = \sin(x+\sin(t))\exp(\cos(y-z))$ so the source function $f$ can be obtained accordingly. We compute both the neural network solutions $\bx_\mathcal{N}$ and $u_\mathcal{N}$ up to time $T = 2$ sequentially by dividing the time interval $[0,2]$ into 10 uniform subintervals so overall 10 steps of time integration are needed to reach the terminal time.

\begin{figure}[h]
\centering
\includegraphics[scale=0.45]{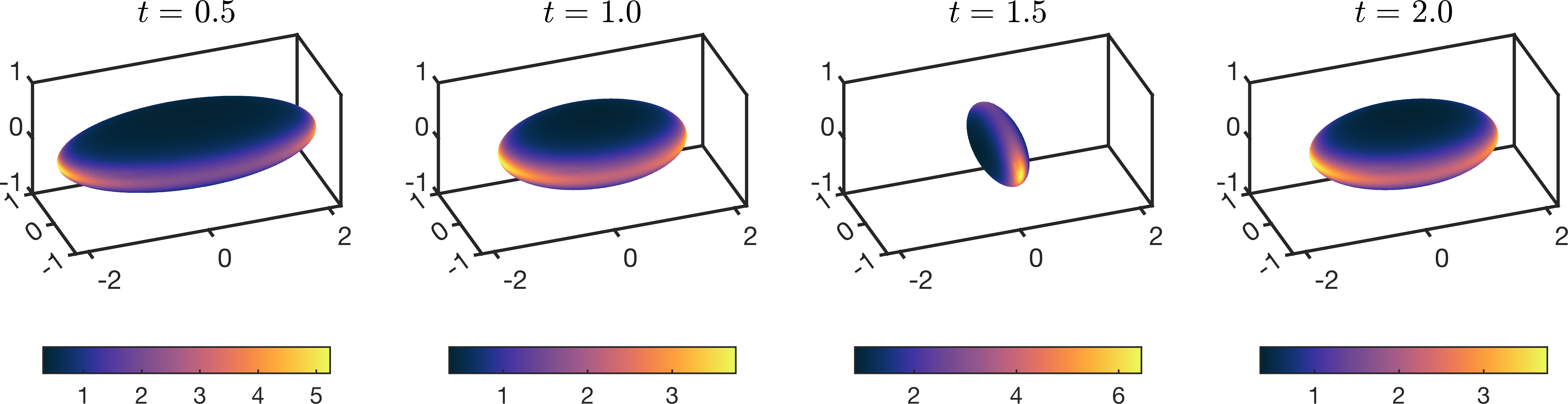}
\caption{The snapshots of the network solution $\bx_\mathcal{N}$ with $N = 40$ neurons at different times. The color code indicates the magnitude of the mean curvature $H_\mathcal{N}$.}
\label{Fig:oscillating_ellipsoid}
\end{figure}

%%%%%%%%%%%%%%%%%%%%%%%%%%%%%%%%%%%%%%%%%
\paragraph{\textbf{Predictive accuracy for tracking the surface}}
%%%%%%%%%%%%%%%%%%%%%%%%%%%%%%%%%%%%%%%%%
The surface evolutional differential equation~(\ref{Eq:surface_evolve}) is solved using the loss model~(\ref{Eq:Loss_surface}).
Given the number of training points $M = 800$ and $M_0 = 100$, we train the network model with different widths $N$ of the hidden layer.
Table~\ref{Table:evolving_surface} reports the relative $L^2$ errors of the network predictive surface configuration $\bx_\mathcal{N}$, the normal vector $\mathbf{n}_\mathcal{N}$, and the mean curvature $H_\mathcal{N}$ at $T=2$.  One can see that the network solution ~(\ref{Eq:NN_surface})  not only predicts accurately for the surface configuration, but also for the normal vector and mean curvature at the test points. Those relative $L^2$ errors range from $O(10^{-4})-O(10^{-6})$ using merely $10-40$ neurons in the hidden layer.
In addition, we depict the snapshots of the predictive surface configuration $\bx_\mathcal{N}$ and mean curvature $H_\mathcal{N}$ for $N = 40$ in Fig.~\ref{Fig:oscillating_ellipsoid}. We should also point out that the present method is mesh-free and the implementation is much easier in comparison with the traditional grid based methods~\cite{HCLT19}.

\begin{table}[h]
\centering
\begin{tabular}{c|ccc}
\hline
$(N,N_p)$    & $\|\bx_\mathcal{N}-\bx\|_2/\|\bx\|_2$ & $\|\mathbf{n}_\mathcal{N}-\mathbf{n}\|_2/\|\mathbf{n}\|_2$ & $\|H_\mathcal{N}-H\|_2/\|H\|_2$ \\
\hline
$(10,80)$     &  4.574E$-$04  &  3.317E$-$04 &  7.360E$-$04\\
$(20,160)$   &  8.678E$-$05  &  3.305E$-$05 &  1.037E$-$04 \\
$(30,240)$   &  4.206E$-$06  &  1.980E$-$06 &  7.156E$-$06 \\
$(40,320)$   &  1.116E$-$06  &  1.892E$-$06 &  7.393E$-$06   \\
\hline
\end{tabular}
\caption{The average relative $L^2$ errors for the surface configuration $\bx$, normal vector $\mathbf{n}$, and mean curvature $H$ at $T=2$. For each case the number of training points is fixed by $M = 800$ and $M_0 = 100$. The total number of learnable parameters for the network expression~(\ref{Eq:NN_surface}) is $N_p = 8N$.}
\label{Table:evolving_surface}
\end{table}

%%%%%%%%%%%%%%%%%%%%%%%%%%%%%%%%%%%%%%%%%
\paragraph{\textbf{Predictive accuracy for solving the advection-diffusion equation on an evolving surface}}
%%%%%%%%%%%%%%%%%%%%%%%%%%%%%%%%%%%%%%%%%
The loss model~(\ref{Eq:Loss_PDE_evolve}) is used to find the solution to the advection-diffusion equation~(\ref{Eq:advection_diffusion}).
Again, using $M = 800$ and $M_0 = 100$ training points, we first train the network with $N = 40$ neurons to predict $\bx_\mathcal{N}$, $\mathbf{n}_\mathcal{N}$, and $H_\mathcal{N}$ at the training points (i.e., at $(\theta^i,\phi^i,t^i)$, we set $\bx^i_\mathcal{N} = \bx^i$, $\mathbf{n}^i_\mathcal{N} = \mathbf{n}(\bx^i,t^i)$, and $H^i_\mathcal{N} = H(\bx^i,t^i)$), and then use them as the inputs in the loss function
(\ref{Eq:Loss_PDE_evolve}).
Table~\ref{Table:PDE_evolving_surface} shows that, simply using the shallow network representation, our proposed solver can indeed achieve high accurate predictions for different number of neurons $N = 10, 20, 30, 40$ used in the hidden layer. We must emphasize that, the present method enjoys the advantage of mesh-free, so the implementation can be easily done to find the embedded solution $u$.

\begin{table}[h]
\centering
\begin{tabular}{cc}
\hline
$(N,N_p)$    &  $\|u_\mathcal{N}-u\|_2/\|u\|_2$\\
\hline
$(10,60)$     &  1.968E$-$03 \\
$(20,120)$   &  6.557E$-$04 \\
$(30,180)$   &  3.903E$-$05 \\
$(40,240)$   &  2.885E$-$05 \\
\hline
\end{tabular}
\caption{The average relative $L^2$ errors for the solution $u$ at $T=2$. For each case the number of training points is fixed by $M = 800$ and $M_b = 100$. The total number of learnable parameters for the network expression is $N_p = 6N$.}
\label{Table:PDE_evolving_surface}
\end{table}

%%%%%%%%%%%%%%%%%%%%%%%%%%%%%%%%%%%%%%%%%
\subsection{Surfactant transport on the droplet surface under shear flow}
\begin{figure}[h]
\centering
\includegraphics[scale=0.47]{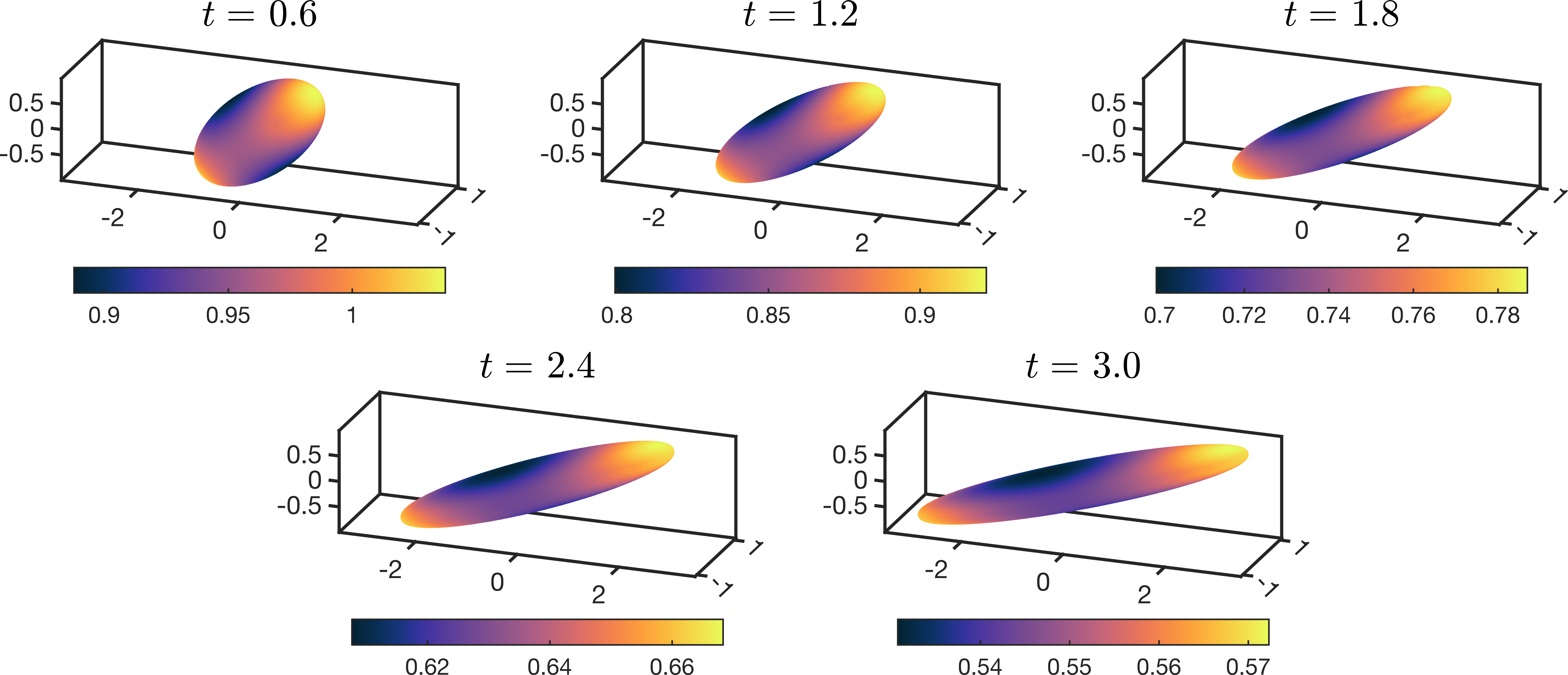}
\caption{The snapshots of the droplet surface configuration $\bx_\mathcal{N}$ ($N = 50$) and surfactant concentration $u_\mathcal{N}$ ($N = 100$) at different times. The color code indicates the magnitude of $u_\mathcal{N}$.}
\label{Fig:shear}
\end{figure}

%%%%%%%%%%%%%%%%%%%%%%%%%%%%%%%%%%%%%%%%%
As an application, we mimic the simulation of surfactant transport on a droplet surface \cite{SHL18, HCLT19} that has been extensively studied using various numerical methods in literature. Here we neglect the fluid effect but simply apply the known shear flow $\mathbf{v}(x,y,z,t) = (z,0,0)$ to the droplet.  The initial shape of the droplet surface is set as a unit sphere located at the origin, and will be elongated by the shear flow along the $x$-direction. One can simply derive the exact surface configuration under this flow as
\begin{align*}
\Gamma(t) = \left\{(x,y,z)\Bigg| \left(x-tz\right)^2 + y^2 + z^2 = 1 \right\}.
\end{align*}
The initial surfactant concentration $u$ is set to be uniform as $u(x,y,z,0) = 1$ while the source term is $f(x,y,z,t) = 0$. We construct the network representation for $\bx_\mathcal{N}$ and $u_\mathcal{N}$ with $N = 50$ and $N = 100$ neurons respectively, and use $M= 1000$ and $M_0 = 500$ training points in the loss models.
The simulation is performed up to time $T = 3$ sequentially by dividing the time interval $[0,3]$ into 10 uniform subintervals so overall 10 steps of time integration are needed to reach the terminal time.
The snapshots for the droplet configuration $\bx_\mathcal{N}$ along with the surfactant concentration $u_\mathcal{N}$ are shown in Fig.~\ref{Fig:shear}. As seen, due to the presence of the applied shear flow, the surfactant is swept toward the both tips of the droplet surface as time evolves, leading to high concentration at the tips while low concentration at the sides of the surface.
This concentration distribution is commonly observed in the presence of shear flow even with the fluid effect ~\cite{HCLT19}.
We should point out that, the present mesh-free neural network method has no difficulty to handle the scenario of large surface distortion (see $t = 3$ in the figure), while in traditional numerical methods, the droplet surface must be re-meshed from time to time to keep accurate and stable computations.

\begin{figure}[h]
\centering
\includegraphics[scale=0.4]{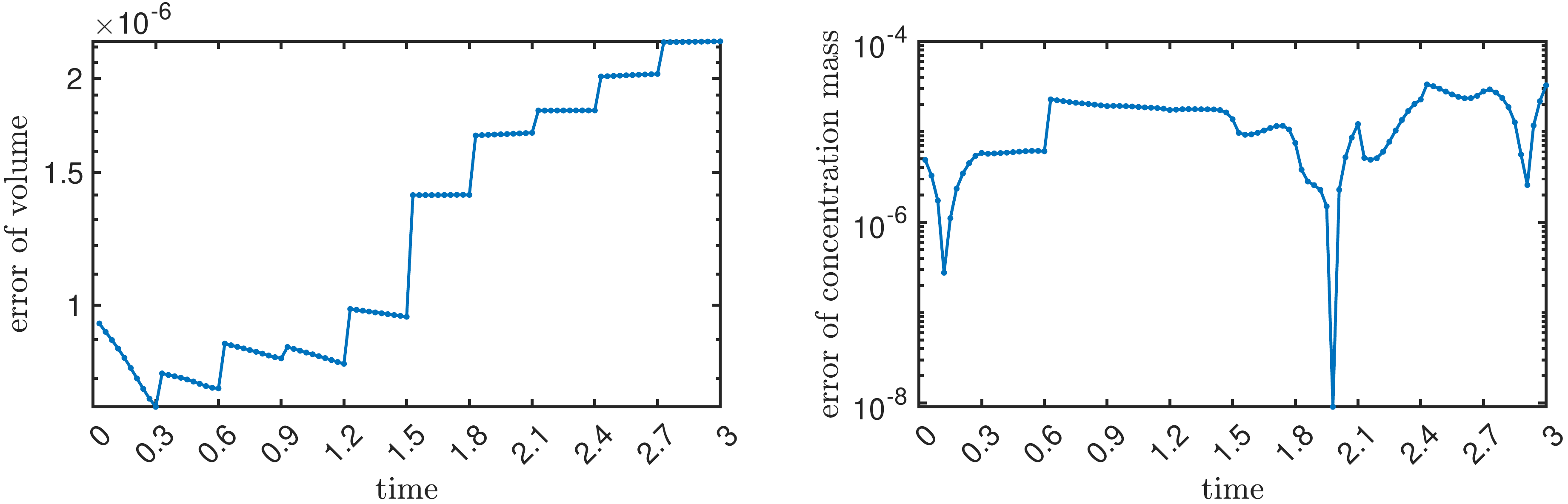}
\caption{The time plots of droplet volume error $|V(t)-V(0)|/V(0)$ (left), and total surfactant mass error $|m(t)-m(0)|/m(0)$ (right).}
\label{Fig:shear_error}
\end{figure}

Since the given flow $\mathbf{v}$ is incompressible, the droplet volume $V(t) = \frac{1}{3}\int_{\Gamma(t)}\bx\cdot\mathbf{n}\mbox{ d}S$ should be conserved as a constant $V(0) = 4\pi/3$.  Meanwhile, without the additional source ($f = 0$), the total surfactant mass $m(t) = \int_{\Gamma(t)} u\mbox{ d}S$ is also conserved as its initial value $m(0) = 4\pi$.
Despite the present method dose not guarantee the numerical conservation for these two quantities, we plot the evolutions of the relative error for droplet volume $|V(t)-V(0)|/V(0)$ and total surfactant mass $|m(t)-m(0)|/m(0)$ in Fig.~\ref{Fig:shear_error}. These two surface integrations for $V(t)$ and $m(t)$ are performed by Gauss-Legendre quadrature rule in $\theta$-direction and midpoint rule in $\phi$-direction.
One can see that both error plots are discontinuous at the endpoint of each time subinterval since the initial conditions for $\bx_\mathcal{N}$ and $u_\mathcal{N}$ are resumed in our loss models.
The relative volume error reaches as low as $O(10^{-6})$ even when the surface is highly distorted at $t = 3$, and the total surfactant mass error reaches $O(10^{-5})$. This results outperform the ones obtained in \cite{HCLT19}.

%%%%%%%%%%%%%%%%%%%%%%%%%%%%%%%%%%%%%%%%%
\section{Conclusion and future works}\label{Sec:conclusion}
%%%%%%%%%%%%%%%%%%%%%%%%%%%%%%%%%%%%%%%%%

In this paper, a completely shallow physics-informed neural network is developed to solve Laplace-Beltrami and diffusion equations on static surfaces, and advection-diffusion equation on evolving surfaces.
Those surface PDEs are written in Eulerian coordinates in which geometrical differentiations are calculated by conventional differential operators. For the static surface case, with the aid of level set function, the surface geometrical quantities such as the normal and mean curvature of the surface can be computed directly and used in our surface differential expressions. The loss function hereby penalizes the equation residual written in the form of Cartesian differential operators instead of imposing normal extension constraints used in literature.  As for the evolving surface, we additionally introduce a prescribed hidden layer to enforce the topological structure of the surface and use the network to learn the homeomorphism between the surface and the prescribed topology.
The proposed network structure is designed to track the surface and solve the equation simultaneously. Since the present neural network uses only one hidden layer, the model is easy to implement and train. Numerical results show high predictive accuracy using just a moderate number of neurons in the hidden layer.

\textcolor{black}{The traditional mesh-free method represents the solution by a linear combination of some chosen radial basis functions (RBFs; for instance, Gaussian), and enforces the solution to satisfy the PDE directly at some chosen points. In fact, one can regard the present shallow neural network solution as a linear combination of activation basis in which the weights and bias must be determined via learning. It would be nice to make a fair performance comparison (including the computational cost and accuracy) between the RBFs method and the present neural network method.  But this is beyond the scope of the paper which we shall leave it as our future work.}
%\textcolor{black}{We take into account the whole Cartesian differential expressions in the loss function. As a result, the proposed neural network approach seeks the embedded solution in the entire $\mathbb{R}^3$ space rather than the parameter space. This step significantly avoids some potential issues caused by the periodicity or pole conditions along the surface.
%Consequently, with the support of universal function expressivity~\cite{Cybenko1989, Hor91, Mhaskar96, RLM21}, using shallow neural network representation to approach smooth solutions in general attains high accuracy predictive results.}

The considered surfaces in this paper are defined with given level set representations. At each training point, both normal vector and mean curvature are required in the present model implementation that can be easily computed by the usage of level set function. As a forthcoming extension, we shall consider PDEs on a point cloud of closed surface in which the level set function is not available. It is apparently a challenging task to compute the normal vectors and mean curvatures at those training points especially when the surface evolves. So it may be worthy to explore a hybrid method that combines machine learning and traditional numerical techniques to tackle the PDEs on evolving surfaces.
%Meanwhile, the considered PDE problems on evolving surfaces in this paper are defined with given prescribed velocity fields. As a forthcoming extension, we shall consider the full dynamic of interfacial flows; namely, the Navier-Stokes equations is taken into account to drive the flow. In addition, we aim to utilize the proposed homeomorphism network representation to solve high-order geometric flow problems such as mean curvature flow or Willmore flow.
We leave this to our future work as well.

%%%%%%%%%%%%%%%%%%%%%%%%%%%%%%%%%%%%%%%%%
\section*{Acknowledgement}
%%%%%%%%%%%%%%%%%%%%%%%%%%%%%%%%%%%%%%%%%

W.-F. Hu, T.-S. Lin and M.-C. Lai acknowledge the supports by National Science and Technology Council, Taiwan, under the research grant 111-2115-M-008-009- MY3, 111-2628-M-A49-008-MY4 and 110-2115-M-A49-011-MY3, respectively.

%%%%%%%%%%%%%%%%%%%%%%%%%%%%%%%%%%%%%%%%%
\section*{Appendix}
%%%%%%%%%%%%%%%%%%%%%%%%%%%%%%%%%%%%%%%%%

Here we present the derivations of the relation between surface differential operators and conventional differential operators in Euclidean space.
We begin by considering the surface gradient operator $\nabla_s$, which describes the changing rate along a regular surface (tangent to the surface) by removing the normal component in conventional gradient
\begin{align*}
\nabla_s u = \nabla u - \partial_n u\,\mathbf{n} = (I - \mathbf{n}\mathbf{n}^T) \nabla u,
\end{align*}
where both $\mathbf{n}$ and $\nabla u$ are aligned as column vectors. On the other hand, the surface divergence operator reads
\begin{align*}
\nabla_s\cdot\mathbf{v} = \left[ (I - \mathbf{n}\mathbf{n}^T) \nabla \right]^T \mathbf{v} = \nabla^T (I - \mathbf{n}\mathbf{n}^T) \mathbf{v} = \nabla\cdot\mathbf{v} - \mathbf{n}^T(\nabla\mathbf{v})\mathbf{n}.
\end{align*}
Combining the above identities, we compute the Laplace-Beltrami operator by
\begin{align*}
\Delta_s u & = \nabla_s\cdot(\nabla_s u) = \nabla_s\cdot(\nabla u - \partial_n u\,\mathbf{n}) \\ &= \nabla\cdot(\nabla u - \partial_n u\,\mathbf{n}) - \mathbf{n}^T\left(\nabla(\nabla u - \partial_n u\,\mathbf{n})\right)\mathbf{n} \\
& = \Delta u - (\nabla\cdot\mathbf{n})\partial_n u - (\nabla \partial_n u)\cdot\mathbf{n} - \mathbf{n}^T(\nabla^2u)\mathbf{n} + \mathbf{n}^T(\nabla(\partial_n u\,\mathbf{n}))\mathbf{n} \\
& = \Delta u - 2H \partial_nu - \mathbf{n}^T(\nabla^2u)\mathbf{n},
\end{align*}
where we have used the fact that $\nabla\cdot\mathbf{n} = 2H$ and $\mathbf{n}^T\nabla\mathbf{n} = \mathbf{0}$.

%%%%%%%%%%%%%%%%%%%%%%%%%%%%%%%%%%%%%%%%%
%\section*{References}
%%%%%%%%%%%%%%%%%%%%%%%%%%%%%%%%%%%%%%%%%

%\bibliography{ref_surface_pde.bib}

\end{document}